\documentclass{amsart}
\usepackage[latin1]{inputenc}
\usepackage{amssymb}
\usepackage{latexsym}
\usepackage{longtable}

\newtheorem{teor}{Theorem}[section]
\newtheorem{cor}[teor]{Corollary}
\newtheorem{lema}[teor]{Lemma}
\newtheorem{prop}[teor]{Proposition}
\newtheorem{conjecture}[teor]{Conjecture}
\newtheorem{defn}[teor]{Definition}
\newtheorem{obswr}[teor]{Observation}
\newtheorem{remarkwr}[teor]{Remark}

\newtheorem{examplewr}[teor]{Example}

\newenvironment{remark}{\begin{remarkwr}\begin{upshape}}{\end{upshape}\end{remarkwr}}
\newenvironment{example}{\begin{examplewr}\begin{upshape}}{\end{upshape}\end{examplewr}}

\newcommand{\Q}{\mathbb{Q}}
\newcommand{\Z}{\mathbb{Z}}

\newcommand{\C}{\mathbb{C}}
\newcommand{\cG}{{\mathbb G}}

\newcommand{\PP}{\mathbb{P}}
\newcommand{\cH}{{\mathbb H}}
\newcommand{\cA}{{\mathbb A}}
\newcommand{\cS}{{\mathbb S}}
\newcommand{\R}{{\mathbb R}}

\newcommand{\M}{\operatorname{M}}
\newcommand{\Gal}{\operatorname{Gal\,}}
\newcommand{\disc}{\mathrm{disc}}
\newcommand{\GL}{\operatorname{GL}}

\newcommand{\End}{\operatorname{End}}

\newcommand{\Res}{\operatorname{Res }}
\newcommand{\NS}{\operatorname{NS}}
\newcommand{\Jac}{\operatorname{Jac}}

\newcommand{\id}{\operatorname{id}}
\newcommand{\Aut}{\operatorname{Aut}}

\newfont{\gotip}{eufb10 at 12pt}

\newcommand{\cO}{{\mathcal O}}

\newcommand{\cL}{{\mathcal L}}
\newcommand{\cF}{{\mathcal F}}
\newcommand{\cR}{{\mathcal R}}
\newcommand{\cQ}{{\mathcal Q}}

\newcommand{\om}{{\omega }}
\newcommand{\ra}{{\rightarrow }}
\newcommand{\qbar}{{\overline{\Q}}}

\newcommand{\PSL}{{\mathrm {PSL}}}
\newcommand{\n}{{\mathrm{n}}}
\newcommand{\tr}{{\mathrm{tr}}}
\newcommand{\Hom}{{\mathrm{Hom }}}

\include{thebibliography}

\begin{document}

\title[Endomorphism algebras of abelian surfaces]
{On finiteness conjectures for endomorphism algebras of abelian
surfaces}

\author{Nils Bruin}
\address{Department of Mathematics, Simon Fraser University, Burnaby, 
BC, Canada V5A~1S6}
\email{bruin@member.ams.org}
\author{E.\ Victor Flynn}
\address{Mathematical Institute, University of Oxford,
Oxford OX1 3LB, United Kingdom} \email{flynn@maths.ox.ac.uk}
\author{Josep Gonz\'{a}lez}
\address{Universitat Polit\`{e}cnica de Catalunya,
Departament de Matem\`{a}tica Aplicada IV (EUPVG), Av.\ Victor
Balaguer s/n, 08800 Vilanova i la Geltr\'{u}, Spain.}
\email{josepg@mat.upc.es}
\author{Victor Rotger}
\address{Universitat Polit\`{e}cnica de Catalunya,
Departament de Matem\`{a}tica Aplicada IV (EUPVG), Av.\ Victor
Balaguer s/n, 08800 Vilanova i la Geltr\'{u}, Spain.}
\email{vrotger@mat.upc.es}

\thanks{The first author is partially supported by an NSERC grant.
The second author is partially supported by EPSRC grant
GR/R82975/01. The third and fourth authors are partially supported
by DGICYT Grant BFM2003-06768-C02-02.}

\begin{abstract}
It is conjectured that there exist only finitely many isomorphism
classes of endomorphism algebras of abelian varieties of bounded
dimension over a number field of bounded degree. We explore this conjecture when
restricted to quaternion endomorphism algebras of abelian
surfaces of $\GL_2$-type over $\Q $ by giving a moduli
interpretation which translates the question into the diophantine
arithmetic of Shimura curves embedded in Hilbert surfaces. We
address the resulting problems on these curves by local and global
methods, including Chabauty techniques on explicit equations of
Shimura curves.
\end{abstract}

\maketitle

\section{Introduction}\label{Conj}
Let $A/K$ be an abelian variety defined over a number field $K$. For
any field extension $L/K$, let $\End _L(A)$ denote the ring of
endomorphisms of $A$ defined over $L$ and $\End^0 _L(A)=\End
_L(A)\otimes \Q $. A conjecture, which may be attributed to Robert
Coleman\footnote{In a personal communication to the last author,
Robert Coleman pointed out that this conjecture was posed by him
during a lecture in a slightly weaker form.}, asserts the following.

$\\ $ {\bf Conjecture} {$\mathbf C(e, g):$} {\em Let $e$, $g\ge 1$
be positive integers. Then, up to isomorphism, there exist only
finitely many rings $\cO $ over $\Z $ such that $\End_L(A)\simeq \cO
$ for some abelian variety $A/K$ of dimension $g$ and a field
extension $L/K$ of a number field $K$ of degree $[K:\Q ]\le e$.}

$\\ $ The conjecture holds in dimension $1$: by the theory of
complex multiplication, if $E/K$ is an elliptic curve and $L/K$ is
an extension of number fields, then $\End_L(E)$ is either $\Z $ or
an order $\cO $ in an imaginary quadratic field $\Q (\sqrt{-d})$
such that the~ring class field $H_{\cO }$ attached to $\cO $ is
contained in $K(\sqrt{-d})$. Thus $[H_{\cO}:\Q ]\leq 2[K:\Q ]$.
Since $[H_{\cO}:\Q ] = 2h(\cO )$, where we let $h(\cO )$ denote the
class number of~$\cO $,~Conjecture C(h, 1) is now a consequence of
the Brauer-Siegel Theorem \cite[Chapter~XVI, Theorem~3]{La}, which
implies that for given $e\ge 1$, there exist finitely many imaginary
quadratic orders $\cO $ such that $h(\cO )\le e$.

Assuming the generalized Riemann hypothesis and using similar ideas,
Greenberg announced \cite{Gr} a generalization of
the above statement to abelian
varieties of arbitrary dimension $g\ge 1$ with multiplication by
orders in complex multiplication fields of degree $2 g$.

Another instance that motivates Coleman's conjecture stems from
the celebrated work of Mazur \cite{Ma}, as we now explain. Let
$E_1$, $E_2/\Q $ be elliptic curves without CM over~$\Q $ and $A=E_1\times
E_2$. Then, it is easily checked that $\End_{\Q }(A) \simeq \Z
\times \Z $ if $E_1$ and $E_2$ are not isogenous over~$\Q $, and
$\End_{\Q }(A) \simeq M_0(N):=\{
\begin{pmatrix} a & b \\ c & d \end{pmatrix}\in \M_2(\Z ), N\mid c \} $
if there is a cyclic isogeny of degree $N$ between $E_1$ and
$E_2$. By \cite[Theorem~1]{Ma} (for $N$ prime) and the discussion
on \cite[p.\,131]{Ma} (for arbitrary $N$), this holds for only
finitely many values of $N$.

As in \cite[p.\,191]{Py}, we say that an abelian variety $A$ defined
over $\Q $ is {\em of $\GL _2$-type over $\Q $} if the endomorphism
algebra $\End _{\Q }^0(A)$ is a number field $E$ of degree $[E:\Q ]
= \dim A$. These abelian varieties had been introduced by Ribet in
\cite[p.\,243]{Ri} (in a slightly more general way) and this
terminology is motivated by the fact that if $E$ is a number field
of degree $\dim A$ which is contained in $\End^0_{\Q }(A)$, then the
action of $\Gal(\qbar /\Q)$ on the $\ell$-adic Tate module
associated with $A$ defines a representation with values in $\GL_2(E
\otimes \Q_{\ell })$. According to \cite[p.\,244]{Ri}, $E$ must be
either a totally real or a complex multiplication number field.

An abelian variety $A$ is called {\em modular over $\Q $} if it is
a quotient of the Jacobian variety $J_1(N)$ of the modular curve
$X_1(N)$ defined over $\Q $. If moreover $A$ is simple over $\Q$,
its modularity over $\Q$ is equivalent to the existence of an
eigenform $f\in S_2(\Gamma_1(N))$ such that $A$ is isogenous over
$\Q$ to the abelian variety $A_f$ attached by Shimura to~$f$. As
is well-known, all simple modular abelian varieties $A$ over $\Q$
are of $\GL _2$-type over $\Q $ and the generalized
Shimura-Taniyama-Weil Conjecture predicts that the converse is
also true (cf.\,e.\,g.\,\cite[p.\,189]{Py}). As was shown by Ribet
in \cite[Theorem 4.4]{Ri}, this conjecture holds if Serre's
Conjecture \cite[Conjecture 3.2.4$_?$]{Se} is assumed.

As we mentioned, Conjecture ${\mathbf C(e, g)}$ is settled for $e\ge
1$, $g=1$. For the particular case $e=1$, we have that if $E/\Q $ is
an elliptic curve over $\Q $ and $L/\Q $ is a field extension, then
$\End_L^0(E)=\Q $ or $\Q (\sqrt {-d})$ for $d=1, 2, 3, 7, 11, 19,
43, 67$ or $163$. On the other hand, the case $g\geq 2$ is
completely open. The aim of this article is to address the question
for quaternion endomorphism algebras of abelian surfaces of
$\GL_2$-type over $\Q $.

In general, it is known that if $A$ is an absolutely simple
abelian surface  of $\GL_2$-type over $\Q$ then, for any number
field $L$, the endomorphism algebra $\End ^0_{L}(A)$ is either  a
real quadratic field or an indefinite division quaternion algebra
over $\Q $ (cf.\,\cite[Proposition 1.1, Theorem 1.2 and
Proposition 1.3]{Py}).

We recall some basic facts on the arithmetic of quaternion
algebras (cf.\,\cite[Ch.\,I $\S 1$ and Ch.\,III, Theorem 3.1]{Vi}
for these and other details). A quaternion algebra $B$ over $\Q $
is a central simple algebra $B$ of rank $4$ over $\Q $. For any
$a, b\in \Q^*$, one can define the quaternion algebra $(\frac{a,
b}{\Q })=\Q +\Q i+\Q j+\Q i j$, where $i^2=a$, $j^2=b$ and $i j =-
j i$. Any quaternion algebra is isomorphic to $(\frac{a, b}{\Q })$
for some $a, b\in \Q ^*$. The reduced discriminant of $B$ is the
square-free integer $D=\mathrm{disc}(B):=\prod p$, where $p$ runs
through the (finitely many) prime numbers such that $B\otimes
\Q_p\not \simeq \M_2(\Q_p)$. Since for any square-free positive
integer $D$ there exists a (single up to isomorphism) quaternion
algebra $B$ with $\disc(B)=D$, we shall denote it by $B_D$. We
have $D=1$ for $B=\mathrm{M}_2(\Q )$, and this is the only non
division quaternion algebra over $\Q $. A quaternion algebra $B$
over $\Q $ is called indefinite if $B\otimes \R \simeq \M_2(\R )$
or, equivalently, if $D$ is the product of an {\em even} number of
prime numbers.

For $\alpha \in B$, let $\bar{\alpha}$ denote its conjugate and
write $\n :B\ra \,\Q $ and $\tr :B\ra \,\Q $ for the reduced norm
and the reduced trace on $B$, respectively. An order $\cO $ in $B$
is a subring of $B$ of rank $4$ over $\Z $ such that $\n(\alpha ),
\tr (\alpha )\in \Z $ for all $\alpha \in \cO $. The order is
maximal if it is not properly contained in any other order. If
$B_D$ is indefinite, there exists a single maximal order in $B_D$
up to conjugation by elements in $B_D^*$, which we will denote by
$\cO_D$.

\begin{defn}
Let $m>1$ and $D=p_1\cdots p_{2r}$ for some $r\ge 1$ be
square-free integers and let $B_D$ be a quaternion algebra over
$\Q $ of discriminant $D$. We say the pair $(D, m)$ is {\em
modular} over $\Q$ if there exists a modular abelian surface $A/\Q
$ such that
$$\End _{\qbar }^0(A)\simeq B_D \mbox{ and } \End _{\Q
}^0(A)\simeq \Q(\sqrt m ).$$

We say the pair $(D, m)$ is {\em premodular} over $\Q$ if there
exists an abelian surface $A$ of $\GL _2$-type over $\Q $ such
that $$\End _{\qbar }^0(A)\simeq B_D \mbox{ and } \End _{\Q
}^0(A)\simeq \Q(\sqrt m ).$$
\end{defn}

We state a particular consequence of Coleman's Conjecture
separately.

\begin{conjecture}\label{QuatConj}
The set of premodular pairs $(D, m)$ over $\Q $ is finite.
\end{conjecture}

It is worth noting that, given a fixed quaternion algebra
$B_D$, there are infinitely many real quadratic fields $\Q(\sqrt m
)$ that embed in $B_D$, since any field $\Q(\sqrt{m})$ with $m$
such that $(\frac{m}{p})\neq 1$ for all $p\vert D$ does embed
in~$B_D$ (cf.\,\cite[Ch.\,III $\S 5$ C]{Vi}). Thus, the finiteness
of premodular pairs $(D, m)$ over $\Q $ for a fixed $D$ is
also not obvious a priori.

A further motivation for Conjecture \ref{QuatConj} is
computational. Define the {\em minimal level} of a modular pair
$(D, m)$ as the minimal $N$ such that there exists a newform $f\in
S_2(\Gamma _0(N))$ with $(B_D,\Q(\sqrt m))\simeq (\End ^0_{\qbar
}(A_f), \End ^0_{\Q }(A_f))$. The computations below are due to
Koike and Hasegawa \cite{Ha} for $N\le 3000$. By means of Steins's
program {\em Hecke} implemented in \cite{magma}, we extended
these computations for $N\leq 7000 $.

\begin{prop}
The only modular pairs $(D, {m})$ of minimal level $N\leq 7000$ are:

$$
\begin{array}{|c|c|c|c|c|c|c|c|}
\hline (D, m) & (6, 2) & (6, 3) & (6, 6) & (10, 10) & (14, 7)  & (15,15) & (22, 11)\\
\hline    N   &   675  &  1568  &  243   & 2700  &     1568  &
3969 & 5408 \\
\hline
\end{array}
$$

\end{prop}
\vskip 0.5cm

In Theorem \ref{MAIN} (iv) we show that the above are not the only
examples of premodular pairs $(D, m)$ over $\Q $. According
to the generalized Shimura-Taniyama-Weil Conjecture in dimension
two, these pairs should actually be {\em modular} pairs.

On the other hand, it is remarkable that not a single example of a
pair $(D, m)$ has ever been excluded from being modular or premodular over $\Q $. In this work we present the first
examples, either obtained by local methods or by methods using
global information, summarised in the following result, which we
shall prove by the end of Section~6.

\begin{teor}\label{MAIN}$\ $
\begin{itemize}
\item[(i)] If $(D, m)$ is a premodular pair over $\Q
$, then $m = D$ or $m=\frac{D}{p}$ for some prime number $p\mid D$
which does not split in $\Q (\sqrt{D/p})$.

\item[(ii)] Let $p, q$ be odd prime numbers. If $(\frac{q}{p})=1$ or
$p\equiv 1$ mod $12$ or $p\equiv q\equiv 1$ mod $4$, then $(p\cdot
q, q)$ is not premodular over $\Q $.

\item[(iii)] The pairs
$$\begin{array}{l}
(D, m)\in \{({10}, 2), (15,3), (15,5),(21,3),(26, 2), ({26},{13}),
(33,11), (38, 2), ({38}, {19}), \\
\hspace{2em} (46,23),(51, 3),(58,2),(91,91),(106,53),(115, 23 ),
(118,59),
(123,123), \\
\hspace{2em} (142,2), (155,
5),(155,31),(155,155),(158,158),(159,3),(202,101),(215, 43),\\
\hspace{2em} (326, 326), (446, 446),(591,3),(1247, 43)\}
\end{array}$$
are not premodular over $\Q $.

\item[(iv)] The pairs
$$\begin{array}{l}
(D,m)\in\{(6,2),(6,3),(6,6),(10,5),(10,10),(14,7),(14,14),(15,15), \\
\hspace{2em}(21,21),(22,2),(22,11),(22,22),(33,33),(34, 34),(46,46),(26,26),\\
\hspace{2em}(38,38),(58, 29),(58, 58)\}
\end{array}
$$
are premodular over $\Q $.

\item[(v)] Let $D>546$. Then there exist only finitely many $\qbar
$-isomorphism classes of abelian surfaces $A$ of $\GL_2$-type over
$\Q $ such that $\End_{\qbar }(A)\simeq \cO_D$.

\item[(vi)] For the pairs
$$\begin{array}{l}
(D,m)\in \{ (6, 2), (6,3), (6,
6), (10, 5), (10, 10), (14,14), (15,15),(21,21),\\
\hspace{2em}(22, 2), (22,11),(22,22),(33, 33),(34,34),(46, 46)\},
\end{array}$$
there exist infinitely many $\qbar $-nonisomorphic abelian
surfaces $A$ defined over $\Q$ such that $\End_{\Q}^0(A)\simeq
\Q(\sqrt m )$ and $\End_{\overline{\Q}}(A)\simeq \cO_D$.

\item[(vii)]  Up to isomorphism over $\qbar $, there exist exactly two
abelian surfaces $A/\Q $ with $\End_{\Q}^0(A)\simeq \Q(\sqrt 7 )$
and $\End_{\overline{\Q}}(A)\simeq \cO_{14}$.
\end{itemize}
\end{teor}

All pairs $(D, m)$ for $D\le 34$ are covered by Theorem
\ref{MAIN}. As a particularly interesting example, we obtain that
there exists no abelian surface $A$ of $\GL_2$-type over $\Q $
such that $\End_{\qbar }^0(A)\simeq B_{155}$: indeed, by Theorem
\ref{MAIN} (i) and (iii) none of the pairs $(155, m)$ are premodular over $\Q $.

Note also that it follows from Theorem \ref{MAIN} (ii) and the
\v{C}ebotarev Density Theorem that there actually exist infinitely
many pairs $(p\cdot q, {q})$ which are not premodular over
$\Q $.

The strategy followed in this paper is to prove that the condition
for a pair $(D,  m )$ to be premodular over $\Q $ is
equivalent to the existence of a point in a suitable subset of the
set of rational points on an Atkin-Lehner quotient of the Shimura
curve canonically attached to $\cO_D$.

The article and the proof of Theorem \ref{MAIN} are organized as
follows: In Section 2 we introduce Shimura curves, Hilbert
surfaces and forgetful maps between them. In Section 3 we use the
diophantine local properties of Shimura curves to prove parts (i)
and (ii) of Theorem \ref{MAIN} as a combination of Theorem
\ref{mod} and Proposition \ref{local}. In Section 4 we prove a
descent result on the field of definition of abelian surfaces with
quaternionic multiplication. In Corollary \ref{546}, we show how
part (v) follows from our results combined with the work in
\cite{Ro}. Finally, in Sections 5 and 6 we prove the remaining
parts of Theorem \ref{MAIN} by means of explicit computations and
Chabauty techniques on explicit equations of Shimura curves.

\section{Towers of Shimura curves and Hilbert surfaces}\label{tow}

We recall some basic facts on Shimura varieties and particularly on
Shimura curves and Hilbert surfaces. Our main references are
\cite[Sections 1 and 2]{Mo}, \cite[Ch.\,III]{BoCa} and
\cite[Sections 1, 2 and 3]{Ed}. Let $\cS = \Res _{\C /\R }(\cG _{m,
\C })$ be the algebraic group over $\R $ obtained by restriction of
scalars of the multiplicative group. A Shimura datum is a pair $(G,
X)$, where $G$ is a connected reductive affine algebraic group over
$\Q $ and $X$ is a $G(\R )$-conjugacy class in the set of morphisms
of algebraic groups $\Hom (\cS , G_{\R })$, as in \cite[Definition
1.4]{Mo}.

Let $\cA _f$ denote the ring of finite adeles of $\Q $. As in
\cite[Section 1.5]{Mo}, for any compact open subgroup $U$ of
$G(\cA _f)$, let
$$
\mathrm{Sh}_U(G,X)(\C ) = G(\Q )\setminus (X \times G(\cA _f))/U,
$$
which has a natural structure of quasi-projective complex
algebraic variety, that we may denote by $\mbox{Sh}_U(G,X)_{\C }$.

Let $(G,X)$ and $(G',X')$ be two Shimura data and let $U$, $U'$ be
compact open subgroups of $G(\cA _f)$ and $G'(\cA _f)$,
respectively. A morphism $f:G\rightarrow G'$ of algebraic groups
which maps $X$ into $X'$ and $U$ into $U'$ induces a morphism
$$\mathrm{Sh}_f: \mathrm{Sh}_{U}(G,X)_{\C }\rightarrow
\mathrm{Sh}_{U'}(G',X')_{\C }$$ of algebraic varieties
(cf.\,\cite[Section 1.6.3]{Mo}).

In this section, we consider two particular instances of Shimura
varieties: Shimura curves attached to an indefinite quaternion
algebra and Hilbert surfaces attached to a real quadratic number
field.

\subsection{Shimura curves}\label{Shimura}

Let $B_D$ be an indefinite quaternion algebra over $\Q $ of
reduced discriminant $D$ and fix an isomorphism $\Phi : B_D\otimes
\R \stackrel{\simeq }{\ra }\mbox{M}_2(\R )$. Let $\cO _D\subset
B_D$ be a maximal order and let $G/\Z $ be the group scheme $\cO
_D^*$. We have that $G(\Q )=B_D^*$ and $G(\cA _f)= \prod _p \cO
_{D, p}^*$, where for any prime number $p$, we let $\cO_{D, p} =
\cO_D \otimes \Z _p$. Let $X=\cH ^{\pm }$ be the $\GL _2(\R
)$-conjugacy class of the map $a+bi\mapsto
\begin{pmatrix}
  a & -b \\
  b & a
\end{pmatrix}$. As complex
analytical spaces, $\cH ^{\pm }$ is the union of two copies of
Poincar\'{e}'s upper half plane $\cH $.

For any compact open subgroup $U$ of $G(\cA _f)$, let $X_{U, \C }
= \mbox{Sh}_U(G, X)_{\C }$ be the Shimura curve attached to the
Shimura datum $(G, X)$ and $U$. It is the union of finitely many
connected components of the form $\Gamma _i\setminus \cH $, where
$\Gamma _i$ are discrete subgroups of $\PSL _2(\R )$.

Fix a choice of an element $\mu \in \cO_D$ such that $\mu
^2+\delta = 0$ for some $\delta \in \Q ^*, \delta > 0$ and let
$\varrho _{\mu }:B_D\ra B_D$, $\beta \mapsto \mu ^{-1}\bar {\beta
}\mu $. For any scheme $S$ over $\C $, let $\cF _{U, \mu }(S)$ be
the set of isomorphism classes of $(A, \iota , \nu , \cL )$, where
$A$ is an abelian scheme over $S$, $\iota : \cO _D\hookrightarrow
\End _S(A)$ is a ring monomorphism, $\nu $ is an $U$-level
structure on $A$ and $\cL $ is a polarization on $A$ such that the
Rosati involution $*:\End^0 _S(A) \rightarrow \End^0_S(A)$ is
$\varrho _{\mu }$ on $B_D$ (cf.\,\cite[p.\,128]{BoCa}). As is well
known, $X_{U, \C }$ coarsely represents the moduli functor $\cF
_{U, \mu }$.

A point $[A, \iota , \nu , \cL ]\in X_{U, \C }(\C )$ is called a
{\em Heegner point} or a {\em CM point} if $\iota $ is not
surjective or, equivalently, if $A$ is isogenous to the square of
an elliptic curve with complex multiplication
(cf.\,\cite[pp.\,16-17]{JoPh}, \cite[Definition 4.3]{Ro2}).

The modular interpretation implies that the reflex field of the
Shimura datum $(G, X)$ is $\Q $ and that $X_{U, \C }$ admits a
canonical model $X_{U, \Q }$ over $\Q $, which is the coarse
moduli space for any of the above moduli functors $\cF _{U, \mu }$
extended to arbitrary bases over $\Q $ (cf.\ \cite[Ch.\,III,
1.1-1.4]{BoCa}, \cite[Section~2]{Mo}). The isomorphism class of
the algebraic curve $X_{U, \Q }$ does not depend on the choice of
$\mu \in \cO_D$, although its moduli interpretation does
depend on~$\mu$. This is
due to the following remarkable property: given a triplet $(A,
\iota , \nu )$ as above, each choice of an element $\mu \in
\cO_D$, $\mu ^2+\delta = 0$, $\delta
>0$, determines a single reduced polarization $\cL_{\mu }$
compatible with $(A, \iota , \nu )$. Given $\mu_1, \mu_2\in
\cO_D$, $\mu_i^2+\delta_i=0$, $\delta_i>0$ for $i=1,2$, the
natural isomorphism between $\cF_{U,\mu_1}$ and $\cF_{U,\mu_2}$ is
provided by the map $[A, \iota , \nu , \cL_{\mu_1}]\mapsto [A,
\iota , \nu , \cL_{\mu_2}]$.

As a particular case, let $\cO \subset \cO _D$ be an integral
order contained in $\cO _D$, and let $\hat {\cO }^* = \prod _p \cO
_p^*$. Let us simply denote by $X_{\cO , \Q }$ the Shimura curve
$X_{U, \Q }$ for $U=\hat {\cO }^*$. Again, for any fixed $\mu \in
B_D^*$, $\mu ^2+\delta = 0$, it admits the following alternative
modular interpretation: $X_{\cO , \Q }$ coarsely represents the
functor $\hat {\cF }_{\cO , \mu }: \mbox{ Sch}/\Q \ra \mbox{ Sets
}$, sending a scheme $S$ over~$\Q $ to the set of isomorphism
classes of triplets $(A, \iota , \cL )$, where $(A, \cL )$ is a
polarized abelian scheme over $S$ as above and $\iota :
B_D\hookrightarrow \End^0_S(A)$ is a ring monomorphism such that
$\iota (B_D)\cap \End _S(A) = \iota (\cO )$.

The Atkin-Lehner group of $X_{\cO ,\Q }$ is the normalizer $W_{\cO
}= \mathrm{Norm}_{B_D^*}(\cO )/\Q^*\cdot \cO^*$. There is a
natural action of $W_{\cO }$ on the functor $\hat {\cF }_{\cO ,
\mu }$: for any $[\om ]\in W_{\cO }$, we have $\om : \hat {\cF
}_{\cO , \mu }(S)\ra \hat {\cF }_{\cO , \mu }(S)$, $(A, \iota ,
\cL )\mapsto (A, \om^{-1}\iota \om , \cL_{\om })$, where $\cL_{\om
}$ denotes the single reduced polarization compatible with $(A,
\om^{-1}\iota \om)$. This action induces a natural immersion
$W_{\cO }\subseteq \Aut _{\Q }(X_{\cO, \Q })$
(cf.\,\cite[Proposition 1.2.6]{JoPh}).

When $\cO_D$ is a maximal order in $B_D$, $W_{\cO_D}\simeq (\Z
/2\Z )^{2 r}$, where $2 r = \# \{ p \mbox{ prime}: p\vert D\}$ is
the number of ramified primes of $D$. A full set of
representatives of $W_{\cO_D}$ is $\{ \om _m: m\vert D, m>0 \}$,
where $\om _m\in \cO _D$, $\n (\om _m)=m$. As elements of
$W_{\cO_D}$, these satisfy $\om _m^2=1$ and $\om _m\cdot \om _n =
\om _{m n}$ for any two coprime divisors $m$, $n\vert D$
(cf.\,\cite[Proposition 1.2.4]{JoPh}, \cite[Section 1]{Ro}).

\subsection{Hilbert surfaces}

Let $F$ be a real quadratic extension of $\Q $, let $R_F$ be its
ring of integers and let $G$ be the $\Z $-group scheme
$\mbox{Res}_{R_F/\Z } (\GL _2(R_F))$. Since $F\otimes \R \simeq \R
^2$, we have $G(\R ) = \GL _2(\R )\times \GL _2(\R )$. Let $X=\cH
^{\pm }\times \cH ^{\pm }$. For any compact open subgroup $U$ of
$G(\cA _f)$, let $H_{U, \C } = \mbox{Sh}_U(G, X)_{\C }$ be the
Hilbert surface attached to the Shimura datum $(G, X)$ and $U$ as in
the first paragraph of \cite[Section 2]{Ed}.

The Hilbert surface $H_{U, \C }$ admits, in the same way as
$X_{U,\C}$, a canonical model $H_{U, \Q }$ over $\Q $ which is the
coarse moduli space of abelian surfaces $(A, j, \nu, \cL )$ together
with a ring homomorphism $j:R_F\ra \End (A)$, a $U$-level structure
and a polarization $\cL $ on $A$ such that $*_{|j(R_F)}$ is the
identity map. When $U$ is the restriction of scalars of $\GL _2(\hat
R)$ for a given quadratic order $R\subseteq R_F$, we write $H_{R, \Q
}:=\mbox{Sh}_U(G, X)_{\Q }$. As in the Shimura curve case, $H_{R, \Q
}$ can also be regarded as the coarse moduli space of polarized
abelian surfaces with real multiplication by $R$ and no level
structure (cf.\,\cite[3.1]{Ed}).

A point $P\in H_{U, \C }(\C )$ is called a {\em Heegner point} or a
{\em CM point} if the underlying abelian surface $A$ has complex
multiplication in the sense of Shimura-Taniyama
(cf.\,\cite[Definition 1.2 and Lemma 6.1]{Ed}): the endomorphism
algebra $\End^0_{\C }(A)$ contains a quartic CM-field.

\subsection{Forgetful maps}

We consider various forgetful maps between Shimura curves and
Hilbert surfaces with level structure.

For any integral quaternion order $\cO $ of $B_D$, let $\hat {\cO
}^*\subseteq \hat {\cO} _D^*$ be the natural inclusion of compact
groups. The identity map on the Shimura data $(\cO _D^*, \cH ^{\pm
})$ induces a morphism $$X_{\cO , \Q }\longrightarrow X_{\cO _D,
\Q }$$ which can be interpreted in terms of moduli as forgetting
the level structure: $[A, \iota , \nu, \cL ]\mapsto [A, \iota ,
\cL ]$.

Similarly, for any quadratic order $R$ of $F$, there is a natural
morphism $$H_{R, \Q }\longrightarrow H_{R_F, \Q }.$$

Finally, let $R\subset \cO $ be a real quadratic order
\emph{optimally embedded} in $\cO $, which means that $R=F\cap \cO
$, and fix an element $\mu \in B_D^*$, $\mu ^2+\delta =0$, $\delta
\in \Q ^*$, $\delta >0$ symmetric with respect to $R$ (that is,
$\varrho _{\mu }|R = 1_R$). Regard $X_{\cO,\Q }$ as representing
the moduli functor $\hat {\cF }_{\cO , \mu }$. Attached to the
pair $(R, \mu )$ there is a distinguished forgetful morphism
$$\begin{matrix} \pi
_{(R, \mu )}: & X_{\cO , \Q } & \longrightarrow & H_{R, \Q }\\
& [A, \iota :\cO \ra \End (A), \cL ] & \mapsto & [A, \iota
_{|R}:R\ra \End (A), \cL ]\end{matrix}$$ of Shimura varieties
which consists on forgetting the ring endomorphism structure in
the moduli interpretation of these varieties.

Let $R'$ be a quadratic order of $F$
optimally embedded in $\cO _D$. Writing $R=R'\cap \cO $, we obtain the
following commutative diagram.
$$
\begin{matrix}
X_{\cO , \Q } & \longrightarrow  & X_{\cO _D,
\Q } \\
\pi _{(R, \mu )}  \downarrow &  & \downarrow \pi _{(R', \mu )}\\
 H_{R, \Q }& \longrightarrow   & H_{R', \Q }
\end{matrix}$$

The main consequence we wish to derive from the above is simply a
translation into terms of moduli of the problem posed in
Section \ref{Conj}.

\begin{prop}\label{hil}

Let $B_D$ be an indefinite division quaternion algebra over $\Q $,
let $\cO _D$ be a maximal order and let $F=\Q (\sqrt{m})$  for
some square-free integer $m>1$.

Assume that, for any order $R$ of $F$, optimally embedded in $\cO
_D$, and $\mu \in B_D^*$ symmetric with respect to $R$, the set of
rational points of $\pi _{(R, \mu )}(X_{\cO _D, \Q })$ in the
Hilbert surface $H_{R, \Q }$ consists entirely of Heegner points.
Then, the pair $(D, m)$ is not premodular over $\Q $.

\end{prop}

{\em Proof.} Let $A/\Q $ be an abelian surface such that
$\End^0_{\Q }(A)= F =\Q (\sqrt{m})$ and $\End^0_{\qbar }(A)= B_D$.
Let $R=\End _{\Q }(A)$ and $\cO = \End _{\qbar }(A)$, which we
will regard as an order in $F$ and an order in $B_D$ respectively.
By construction, the order $R$ is optimally embedded in $\cO $.
Since $A$ is projective over $\Q $, it admits a (possibly
non-principal) polarization $\cL $ defined over $\Q $. Let $*$
denote the Rosati involution on $B_D$ induced by $\cL $. By
\cite[Theorem 1.2 (4)]{Ro3}, we have $* = \varrho _{\mu }$ for
some $\mu \in B_D^*$ with $\mu ^2+\delta =0$ for some $\delta \in
\Q ^*$, $\delta
>0$. By choosing an explicit isomorphism $\iota:\cO
\stackrel{\sim}{\longrightarrow } \End _{\qbar }(A)$, the triplet
$(A, \iota , \cL )$ produces a point $P$ in $X_{\cO , \Q }(\qbar
)$, when we regard the Shimura curve as coarsely representing the
functor $\hat {\cF }_{\cO , \mu }$.

Moreover, we have $\iota _{|R}:R\simeq \End _{\Q }(A)$. From the
fact that $\cL$ is defined over $\Q $, it follows that $*_{|R}$ is
an anti-involution on $R$. Since $R$ is totally real, it follows
that $*_{|R}$ is the identity. Hence, the point $P$ is mapped to a
point $P_R\in H_{R, \Q }(\Q )$ by the forgetful map $\pi _{(R, \mu
)}:X_{\cO , \Q }\ra H_{R, \Q }$.

Let $\cO _D$ be a maximal order in $B_D$ containing $\cO $ and let
$R'= F\cap \cO _D$, where we regard $F = R\otimes \Q $ as
naturally embedded in $B_D$. By the above commutative diagram of
Shimura varieties, we obtain a point $P_{R'}\in H_{R', \Q }(\Q )$
which lies in the image of the forgetful map $X_{\cO _D, \Q }\ra
H_{R', \Q }$.

Since $\End^0_{\qbar }(A)\simeq B_D$ is a quaternion algebra, it
contains no quartic CM-fields and thus $P_{R'}$ is not a Heegner
point. This proves the proposition. $\Box $

$\\ $ The relevance of Proposition \ref{hil} to our problem is the
following. Firstly, it translates the condition for a pair $(D,
m)$ to be premodular over $\Q $ into the existence of a
suitable rational point on a projection of a Shimura curve.
Secondly, note that $(D, m)$ is a premodular pair over $\Q $
if there exists an abelian surface $A/\Q $ such that $\End_{\Q
}(A)$ is an order in $\Q (\sqrt{m})$ and $\End_{\qbar }(A)$ is an
order in the quaternion algebra $B_D$. Proposition \ref{hil}
reduces our problem to study the set of rational points on the
Shimura curve $X_{\cO_D, \Q }$ attached to a {\em maximal} order
in $B_D$. These curves have been extensively studied, rather than
the more general curves $X_{\cO }$ attached to an arbitrary
quaternion order.

\section{Atkin-Lehner quotients of Shimura curves}\label{Atkin}

Fix a maximal order $\cO _D$ in an indefinite division quaternion
algebra $B_D$ of discriminant $D$ and let us simply denote $X_D =
X_{\cO _D, \Q }$. As explained in Section \ref{Shimura}, $X_D$ is
equipped with the Atkin-Lehner group of involutions $W_D=\{ \om_m :
m\mid D \}\subseteq \Aut _{\Q }(X_D)$. For $m\mid D$, let
$X_D^{(m)}$ be the quotient curve $X_D/\langle \om _m\rangle $ and
$\pi_m:X_D\ra X_D^{(m)}$ the natural projection map.

For any extension field $K/\Q $, let $X_D(K)_h$ denote the subset
of Heegner points of $X_D(K)$ and let $X_D(K)_{nh} =
X_D(K)\backslash X_D(K)_h$ the set of non-Heegner points over $K$.
Similarly, set $X_D^{(m)}(K)_h = \pi_m(X_D(\qbar )_h)\cap
X_D^{(m)}(K)$ and $X_D^{(m)}(K)_{nh} = X_D^{(m)}(K)\backslash
X_D^{(m)}(K)_h$.

\begin{prop}\label{Thm0}

Let $X_D$ be the Shimura curve of discriminant $D$ attached to the
maximal order $\cO_D$ as above. Then,

\begin{enumerate}
\item[(i)]$X_D(\R )=\emptyset $.

\item[(ii)] There exists no abelian surface $A/\R $ such that $\End_{\R
}(A)\supseteq \cO_D$.

\end{enumerate}
\end{prop}

{\em Proof.} (i) is \cite[Theorem 0]{Sh2} when particularized to
Shimura curves. (ii) follows from (i) and the moduli
interpretation of $X_D$ described in Section \ref{Shimura}. $\Box
$

\begin{teor}\label{mod}

Let $m>1$ be a  square-free integer. Assume that the pair $(D, m)$
is premodular over $\Q $. Then,

\begin{enumerate}
\item[(i)] $m\vert D$ and all prime divisors $p\vert \frac{D}{m}$
do not split in $\Q(\sqrt m )$.

\item[(ii)] $X_D^{(m)}(\Q )_{nh}\ne \emptyset $.

\end{enumerate}

\end{teor}

{\em Proof.} Assume that $(D, m)$ is premodular over $\Q$.
By Proposition \ref{hil}, there exists an order $R$ of $F=\Q(\sqrt
m)$ optimally embedded in $\cO _D$ and $\mu \in \cO_D$,
$\mu^2+\delta =0$, $\delta >0$, symmetric with respect to $R$ such
that the set of rational points of $\pi _{(R, \mu )}(X_{\cO _D, \Q
})$ in the Hilbert surface $H_{R, \Q }$ contains a non-Heegner
point.

Assume first that $m\nmid D$. It was shown in \cite[Theorem
4.4]{Ro2}, (cf.\,also \cite[Section~6]{Ro3} when $\delta =D$) that
there is then a birational equivalence
$$\pi _{(R, \mu )}(X_{\cO _D, \Q })\stackrel {\sim }{\dashrightarrow }X_D.$$

This birational morphism is defined over $\Q $ and becomes a
regular isomorphism when restricted to the set of non-Heegner
points. By Proposition \ref{Thm0}, we obtain a contradiction.

Assume now that $m\vert D $. Since the anti-involution
$\varrho_{\mu }$ restricts to the identity map on $R$, we have
that $B_D\simeq (\frac{-\delta ,m}{\Q })$. Again, it follows from
\cite[Theorem 4.4]{Ro2}, that there is a birational equivalence
$$\pi _{(R, \mu )}(X_{\cO _D, \Q })\stackrel
{\sim }{\dashrightarrow }X_D^{(m)} $$ which is defined over $\Q $
and that becomes a regular isomorphism when restricted to the set
of non-Heegner points. We conclude that $X_D^{(m)} (\Q )$ must
contain a non-Heegner point. Moreover, since $F$ must embed in
$B$, \cite[Ch.\,III $\S 5$ C]{Vi} applies to ensure that all
primes $p\vert D$ do not split in $F$. $\Box $

As an immediate consequence of (i), we obtain the following corollary.

\begin{cor}\label{cc}
Given a discriminant $D$ of a division quaternion algebra over
$\Q$, the set of modular pairs $(D, m)$ is finite.
\end{cor}

In view of Theorem \ref{mod}, the diophantine properties of these
curves are crucial for the understanding of Conjecture
\ref{QuatConj}. We first study under what circumstances the curve
$X_D^{(m)}$ has no points over some completion of $\Q $.

Proposition below is \cite[Theorem 2.7]{RSY}. Part (i) has also
been shown in \cite{Cl} by using supersingular abelian surfaces.

\begin{prop}\label{local}
Let $X^{(m)}_D$ be as above. Then $X_D^{(m)}(\Q_v)\ne \emptyset $
for all places $v$ of $\Q $ if and only if one of the following
conditions holds:

\begin{enumerate}

\item[$(i)$] $m=D$

\item[$(ii)$] $m = D/\ell $ for a prime $\ell \not = 2$ such that

\begin{itemize}

\item $(\frac {m}{\ell }) = -1$;

\item {\rm (a)} $(\frac {-m}{\ell })=1$,
$(\frac {-\ell\, }{p} )\not =1$ for all primes $p\mid m $, or {\rm
(b)} $\ell \equiv 1$ {\rm mod} $4$, $p\not \equiv 1$ {\rm mod} $4$
for all primes $p\mid m $;

\item if $r\ge 2$, then we have
$(\frac{-m/p}{\ell })=-1$ for all odd primes $p\mid m $, and if
$2|D$ we also have either $(\frac{-m/2}{\ell })=-1$, or $q\equiv
3$ {\rm mod} $4$ for all primes $q\mid D/2$;

\item for every prime $p\nmid D$, $p<4g^2$, there exists some
imaginary quadratic field that splits $B_D$ and contains an
integral element of norm $p$ or $p m$.

\end{itemize}

\item[$(iii)$] $m = D/2$ such that

\begin{itemize}

\item $m\not\equiv 1$ {\rm mod} $8$;

\item
$p\equiv 3$ {\rm mod} $4$ for all $p\mid m$, or $p\equiv 5$ or $7$
{\rm mod} $8$ for all $p\mid m$;

\item if $r\ge 2$, then for every prime $p\mid m$ we have
$m/p\not\equiv -1$ {\rm mod} $8$;

\item for every prime $p\nmid D$, $p<4g^2$, there exists some
imaginary quadratic field that splits $B_D$ and contains an
integral element of norm $p$ or $p m$.

\end{itemize}

\end{enumerate}

\end{prop}

As a direct consequence of the combination of Theorem \ref{mod}
and Proposition \ref{local}, we obtain Theorem \ref{MAIN} (i).
When $D=p\cdot q$ is the product of two primes, the congruence
conditions of parts (ii) and (iii) above simplify notably, as we
state in Theorem \ref{MAIN} (ii).

As a consequence of part (i), we conclude that global methods alone
may enable us to prove that any pair $(D,{D})$ is not premodular. 
We also note that the last item of parts (ii) and (iii)
of Proposition \ref{local} allows us to produce isolated examples
of pairs like $(159, 3), (215, 43), (591, 3)$ and $(1247, 43)$
which are not premodular over $\Q $ and are not covered by
Theorem \ref{MAIN} (ii).

Finally we remark that \cite{RSY} studies the failure of the Hasse
principle for curves $X_{p q}^{(q)}$ over $\Q $ for suitable
collections of pairs of primes $p, q$. As an example, it is shown
in \cite[Section 3]{RSY} that $X_{23\cdot 107}^{(107)}(\Q
)=\emptyset $ although it does have rational points everywhere
locally. Hence, we obtain that $(23\cdot 107, 107)$ is not premodular over $\Q $ but this can not be derived from Theorem
\ref{MAIN} (ii).

\section{A descent theorem on rational models of abelian surfaces with
quaternionic multiplication}\label{descent}

Let $A/K$ be an abelian variety defined over a number field $K$
and let $\cL $ be a polarization on it. Let $\qbar $ be a fixed
algebraic closure of $\Q $ containing $K$. Let $K_0$ be the field
of moduli of $(A, \cL )$, that is, the minimal subfield $K_0$ of
$K$ such that for each $\sigma \in \Gal(\qbar /K_0)$ there exists
an isomorphism $\mu_{\sigma }: (A^{\sigma }, \cL^{\sigma })\ra
\,(A, \cL )$ of polarized abelian varieties over $\qbar $. A
similar definition works for a triplet $(A, \iota , \cL )$ where
$\iota :R\hookrightarrow \End (A)$ is a monomorphism of rings for
a given ring $R$, by asking the isomorphisms $\mu _{\sigma }$ to
be compatible with the action of $R$ on $A$
(cf.\,\cite[Section~1]{Ro4}, for more details).

The following result is due to Weil \cite[Theorem 3]{We}. See also
the first paragraph of \cite[Section 8]{Ri} for the specific
statement for abelian varieties and a generalization in the
category of abelian varieties up to isogeny.

\begin{prop}\label{Weil} A polarized abelian variety $(A, \cL )/K$ admits a
model over its field of moduli $K_0$ if and only if for each
$\sigma \in \Gal(\qbar /K_0)$ there exists an isomorphism
$\mu_{\sigma }: (A^{\sigma }, \cL^{\sigma })\ra (A, \cL )$ such
that $\mu_{\sigma }\mu_{\tau }^{\sigma } = \mu_{\sigma \tau}$ for
any $\sigma , \tau \in \Gal(\qbar /K_0)$.

\end{prop}

Let $B_D$ be an indefinite division quaternion algebra of reduced
discriminant $D=p_1\cdots p_{2r}$, $r\ge 1$, and let $\cO_D$ be a
maximal order in $B_D$. Let $m\ge 1$, $m\mid D$.

\begin{lema}\label{Aut}
Let $(A, \cL )/K$ be a polarized abelian surface such that
$\End_{\qbar }(A)\simeq \cO_D$. Then $\Aut_{\qbar }(A, \cL )=\{
\pm 1\}$.
\end{lema}

{\em Proof.} Since $\End_{\qbar }(A)\simeq \cO_D$, it follows that
$\Aut_{\qbar }(A)\simeq \cO_D^* = \{ \alpha \in \cO_D: n(\alpha
)=\pm 1\}$, which is an infinite group (cf.\,\cite[Ch.\,IV,
Theorem 1.1]{Vi}). Theorem 2.2 in \cite{Ro1} attaches an element
$\mu (\cL )\in B_D$ to the polarization $\cL $ and shows that
$\Aut_{\qbar }(A, \cL )\simeq \{ \alpha \in \cO_D^*: \bar {\alpha
} \mu (\cL ) \alpha = \mu (\cL )\} $ (cf.\,also \cite[Theorem 1.2
(1)]{Ro3}). By \cite[Theorem 1.2 (2-3)]{Ro3}, $\mu (\cL )^2+d=0$
for some $d\in \Z $, $d\ge D$.

Let $\alpha \in \Aut_{\qbar }(A, \cL )$. We first observe that
$n(\alpha )=1$: Indeed, if $n(\alpha )=-1$ then $\bar {\alpha
}=-\alpha ^{-1}$ and thus $\mu (\cL )\alpha = -\alpha \mu(\cL )$.
This implies that $B_D\simeq (\frac{-1, -d}{\Q })$, which
contradicts the indefiniteness of $B_D$.

Hence, $\Aut_{\qbar }(A, \cL )\simeq \{ \alpha \in \cO_D^*:
n(\alpha )=1, \mu (\cL )\alpha = \alpha \mu(\cL )\} \simeq S^*$,
where $S=\cO_D\cap \Q (\mu(\cL ))$. Since $S$ is an imaginary
quadratic order in $\Q (\mu(\cL ))\simeq \Q(\sqrt{-d})$ for $d\ge
D\ge 6$, we obtain that $\Aut_{\qbar }(A, \cL )=\{ \pm 1\} $.
$\Box $

The next Lemma is Theorem 3.4 (C. (1)) of \cite{DiRo}.

\begin{lema}\label{pp}
Let $A/\Q $ be an abelian surface such that $\End _{K}(A)\simeq
\cO _D$ over an imaginary quadratic field $K$ and $\End _{\Q
}^0(A)\simeq \Q (\sqrt {m})$. Then $A$ admits a polarization $\cL
\in H^0(\Gal (\qbar /\Q ), \NS (A_{\qbar }))$ of degree $d
>0$ if and only if $B_D \simeq (\frac {-D d, m}{\Q })$.
\end{lema}

The next Lemma is essentially due to Ribet.

\begin{lema}\label{Cre}
Let $A/\Q $ be an abelian surface such that $\End^0_{\Q}(A)\simeq
\Q(\sqrt{m})$ and $\End^0_K(A)\simeq B_D$, where
$K=\Q(\sqrt{-\delta })$, $\delta >0$. Then $B_D\simeq
(\frac{-\delta,m}{\Q })$.
\end{lema}

{\em Proof.} This is stated verbatim in \cite[Theorem 1]{Cr}.
However, note that the statement of \cite[Theorem 1]{Cr} is
restricted to {\em modular} abelian surfaces. By applying and
making explicit \cite[Theorem 5.6]{Ri}, it is shown in
\cite[p.\,133]{GLQ}, that the same formula still holds for
arbitrary abelian surfaces of $\GL_2$-type over $\Q $. $\Box $

In \cite[Theorem 4.2]{Mu}, Murabayashi proved a descent result for
principally polarized simple abelian surfaces with quaternionic
multiplication under certain hypotheses. We give an alternative
proof of his result that allows us to generalize it to arbitrarily
polarized abelian surfaces and which is unconditionally valid.

\begin{teor}\label{teodes}

There exists an abelian surface $A/\Q $ such that $\End^0_{\Q
}(A)\simeq \Q (\sqrt{m})$ and $\End_{\qbar }(A)\simeq \cO_D$ if
and only if there exists $Q\in X_D^{(m)}(\Q )_{nh}$ such that
$\pi_m^{-1}(Q)\subset X_D(K)$ for an imaginary quadratic field
$K=\Q (\sqrt{-\delta })$ with
$$B_D\simeq (\frac{-\delta , m}{\Q }).$$

\end{teor}

{\em Proof.} {\bf I.} Assume first that there exists an abelian
surface $A/\Q $ such that $\End^0_{\Q }(A)\simeq \Q (\sqrt{m})$
and $\End_{\qbar }(A)\simeq \cO_D$.

Let us show how can one attach to $A$ a point $Q\in X_D^{(m)}(\Q
)_{nh}$. Let $d\ge 1$ be the minimal integer such that $B_D\simeq
(\frac{-D d, m}{\Q })$. We know from Lemma \ref{pp} that there
exists a polarization $\cL $ on $A$ of degree $d$ defined over $\Q
$. Fix an isomorphism $\iota : \cO_D\stackrel{\sim }{\ra}\End
_{\qbar }(A)$ and let $R=\iota^{-1}(\End_{\Q }(A))$, which is
isomorphic to a quadratic order of $\Q(\sqrt{m})$. Since
$\End^0_{\Q }(A)$ is a real quadratic field, \cite[Theorem 3.4
{\bf C}]{DiRo} shows that $\End _{\qbar }(A) = \End_K(A)=\cO_D$
for some imaginary quadratic field $K=\Q (\sqrt{-\delta })$.

Let $*:B_D\ra B_D$ be the Rosati involution induced by $\cL $. By
\cite[Theorem 1.2 (2-3-4)]{Ro3}, $\beta ^* = \varrho _{\mu }(\beta
)=\mu ^{-1}\bar {\beta }\mu $ for some $\mu \in \cO _D$, $\mu ^2+D
d=0$. If we regard the Shimura curve $X_D$ as coarsely
representing $\hat {\cF }_{\cO_D, \mu }$, the triplet $[A, \iota ,
\cL ]$ produces a point $P$ in $X_D(K)$. Note however that the
triplet $(A, \iota _{|R}, \cL )$ is defined over $\Q $. Hence
$\pi_{(R, \mu )}(P)\in H_{R, \Q }(\Q )$. Since $\varrho _{\mu }$
is symmetric with respect to $R$, as in the proof of Theorem
\ref{mod}, the forgetful morphism $\pi _{(R, \mu )}:X_D\rightarrow
\,H_{R, \Q }$ is birationally equivalent to the composition of the
projection $\pi_m :X_D\ra X_D^{(m)}$ and an immersion of
$X_D^{(m)}$ into $H_{R, \Q }$. To be more precise, as stated in
\cite[Theorem 3.5]{Ro3}, there exists a possibly singular curve
$\tilde X\stackrel{j}{\hookrightarrow }H_{R, \Q }$ embedded in
$H_{R, \Q }$ such that $\pi _{(R, \mu )} =j\cdot b\cdot \pi _m$
where $b:X_D^{(m)}\ra \tilde X$ is a birational map that becomes
an isomorphism away from the set of Heegner points. Hence
$Q=\pi_m(P)\in X_D^{(m)}(\Q )_{nh}$. Finally, by Lemma \ref{Cre},
$B_D\simeq (\frac{-\delta , m}{\Q })$.

 {\bf II.} Conversely, let $K=\Q
(\sqrt{-\delta })$ be an imaginary quadratic field.  Assume that
$B_D\simeq (\frac {-\delta, m}{\Q })$ and let $P\in X_D(K)$ such
that $Q=\pi_m(P)\in X_D^{(m)}(\Q )_{nh}$ (and thus
$\pi_m^{-1}(Q)=\{ P, \om_m(P)\} \subset X_D(K)$). Choose $\mu ,
\om \in \cO _D$, $\mu ^2=-\delta $, $\om ^2=m$, $\mu \om = -\om
\mu $ and let $R = \Q (\om )\cap \cO _D$. Regard the Shimura curve
$X_D$ as coarsely representing $\hat {\cF }_{\cO_D , \mu }$. The
element $\mu $ determines an embedding of $K$ into $B_D$. As is
stated in \cite[Theorem 2.1.3]{JoPh}, the point $P\in X_D(K)$ can
be represented by the $\qbar $-isomorphism class of a polarized
simple abelian surface with quaternionic multiplication $(A_0,
\iota _0, \cL _0)$ completely defined over $K$ and such that the
Rosati involution that $\cL _0$ induces on $B_D$ is $\varrho _{\mu
}$.

As in Part I, the condition $\pi_m (P)\in X_D^{(m)}(\Q )_{nh}$
implies that $\pi _{(R, \mu )}(P)\in H_{R, \Q }(\Q )$, and this
amounts to saying that the field of moduli of $(A_0, \iota
_{0|_R}, \cL _0)$ is~$\Q $.

For any number field $F\subset \qbar $, let $G_F=\Gal(\qbar /F)$.
Let $\sigma\in G_{\Q}\backslash G_K$. Then there exists an
isomorphism $\nu: A_0 \rightarrow A_0^{\sigma}$ such that
$\nu^{*}(\cL_0^{\sigma})=\cL _0$ and $\nu \cdot \om ^{-1}\cdot
\alpha \cdot \om =\alpha ^{\sigma}\cdot \nu$ for all endomorphisms
$\alpha \in B_D=\End^0_K(A)$. In particular,
$$
\nu \cdot \om = \om ^{\sigma}\cdot \nu \,,\quad \nu \cdot \mu
=-\mu ^{\sigma}\cdot \nu\,.
$$
We split the proof into two  parts.

\medskip
\noindent \emph{Step 1}: We show that $\nu $ may be assumed to be defined over $K$.

To prove this claim, we first note that since $\End_{\qbar
}(A_0)=\cO_D$, $\Aut _K(A_0, \cL _0)=\{ \pm 1\} $ by Lemma
\ref{Aut}. Let $\rho_{\nu}:G_K\rightarrow \Aut _K(A_0, \cL _0)=\{
\pm 1\} $ be the group homomorphism defined by
$\rho_{\nu}(\tau)=\nu^{-1}\cdot \nu^{\tau}$.

Suppose that $\nu $ was not defined over $K$, that is
$\rho_{\nu}(G_K)=\{\pm 1\}$. Let $L/K$ be the quadratic extension
such that $G_L=\ker \rho_{\nu}$. Since $L$ is the minimal field of
definition of all homomorphisms in $\Hom (A_0, A_0^{\sigma})$ and
$\Hom (A_0^{\sigma}, A_0)$, we deduce that $L/\Q$ is a Galois
extension. Since $K$ is imaginary, $L/\Q$ can not be cyclic and
there exists a square-free integer $d>1$ such that $L=K(\sqrt
{d})$.

Let $V_K=H^0(A_0,\Omega^1_{A_0/K})$ denote the vector space of
regular differentials on $A$ over $K$. Since $B_D\simeq
\End_K^0(A)$, the action of the endomorphisms on $V_K$ induces an
embedding $*:B_D\hookrightarrow \End_K(V_K)\simeq \M_2(K)$ and an
isomorphism $B_D\otimes K \simeq K+K\mu +K\nu +K\mu \nu \simeq
\M_2(K)$. We may choose basis of $V_K$ such that the matrix
expressions of $\om ^*$ and $\mu ^*$ acting on $V_K$ are
$$
M_{m}=\left (\begin{array}{cc} 0&1\\m&0
\end{array}\right)\,, \quad M_{\delta }=\left (\begin{array}{cc}
\sqrt{-\delta}&0\\0&-\sqrt{-\delta}
\end{array}\right),
$$
respectively. Indeed, this is possible because, as stated by the
Skolem-Noether Theorem (cf.\,\cite[Ch.\,I, Theorem 2.1]{Vi}), all
automorphisms of $\M_2(K)$ are inner and we have $B_D\otimes K\simeq
\M_2(K) = K+K M_{\delta }+K M_{m} + K M_{\delta } M_{m}$.

Let $N\in \GL_2(K(\sqrt d))$ be the matrix expression of $\nu \in
\Hom (A_0, A_0^{\sigma })$ with respect to this basis of $V_K$ and
its Galois conjugate of $V_K^\sigma $. Then $N$ satisfies
$$
N^{\tau}=-N\,,\quad M_m\cdot N = N\cdot M_m^{\sigma} = N\cdot M_m
\,, \quad M_{\delta }\cdot N = -N\cdot M_{\delta }^{\sigma} =
N\cdot M_{\delta }\,,
$$
for $\tau \in G_K\setminus G_L$. Hence, $N=\sqrt {d} \left
(\begin{array}{cc} \beta&0\\0&\beta
\end{array}\right)\,, \beta\in K\,.$ Fix $\sigma $ in
$G_{\Q(\sqrt d)}$, $\sigma \not \in G_K$. We have
$\nu^{\sigma}\cdot \nu \in \Aut (A_0, \cL _0)=\{ \pm 1 \}$, thus
$N\cdot N^{\sigma}=\pm \id$ and $\beta \cdot \beta^{\sigma}=1/d$.
Hence, the normal closure $F$ of $K(\sqrt{\beta})/\Q$ is dihedral
containing $K(\sqrt d)$ and $F/\Q(\sqrt{-d\cdot \delta})$ is
cyclic. Let $\rho_{\beta}: G_K\rightarrow \{ \pm 1\}$ be the
surjective morphism such that $\ker
\rho_{\beta}=G_{K(\sqrt{\beta})}$. Attached to the cocycle $\rho
_{\beta }\in H^1(G_K, \{ \pm 1\})$ there is a polarized abelian
surface $(A_1, \cL_1)$ defined over $K$ together with an
isomorphism $\lambda: (A_0, \cL_0)\rightarrow (A_1, \cL_1)$ such
that $\lambda^{\tau}=\lambda \cdot \rho_{\beta}(\tau)$. We claim
that $\phi=\lambda^{\sigma}\cdot \nu \cdot
\lambda^{-1}:A_1\rightarrow A_1^{\sigma}$ is defined over $K$.
Indeed, for any $\tau\in G_K$,

$$
\phi^{\tau}=(\lambda^{\sigma \cdot \tau \cdot
\sigma^{-1}})^{\sigma}\cdot \nu^{\tau}\cdot (\lambda^{-1})^{\tau}
= \rho_{\beta{}}(\sigma \cdot \tau \cdot \sigma^{-1}\cdot
\tau^{-1})\cdot \rho_{\nu}(\tau)\cdot \phi\,.
$$
Since $\sigma \cdot \tau \cdot \sigma^{-1} \cdot \tau^{-1}\in
G_{K(\sqrt{\beta})}$ if and only if $\tau \in G_{K(\sqrt{d})}$, we
obtain that $\phi^{\tau}=\phi$.

Moreover, all endomorphisms of $A_1$ are of the form $\lambda
\cdot \varphi \cdot \lambda^{-1}$ with $\varphi $ in $\End_K
(A_0)$. These are all defined over $K$ because $(\lambda \cdot
\varphi \cdot \lambda^{-1})^{\tau}= \delta_{\beta}(\tau \cdot
\tau^{-1})\lambda \cdot \varphi \cdot \lambda^{-1} = \lambda \cdot
\varphi \cdot \lambda^{-1}\,.$

We therefore assume that $\nu $ is defined over $K$.

\medskip
\noindent \emph{Step 2}: We show that $(A_0, \cL_0)$ admits a
model over $\Q $ with all its endomorphisms defined over $K$.

We do so by applying Proposition \ref{Weil}. Since $\nu
^{\sigma}\cdot \nu \in \Aut (A_0, \cL _0)$, we have $\nu
^{\sigma}\cdot \nu = \epsilon \id$ with $\epsilon\in\{\pm 1 \}$.
Using the same basis of $H^0(A_0,\Omega^1_{A_0/K})$ and
$H^0(A_0^{\sigma},\Omega^1_{A_0^{\sigma}/K})$ as above, the matrix
expression $N\in\GL_2(K)$ of $\nu $ is such that $ M_m\cdot N =
N\cdot M_m^{\sigma} = N\cdot M_m \,, \quad M_{\delta }\cdot N
 =-N\cdot M_{\delta }^{\sigma} = N \cdot M_{\delta }\,.$
It follows that
$$
N= \left (\begin{array}{cc} \beta&0\\0&\beta
\end{array}\right)\,,\quad \beta\in K\,.
$$
Hence, $\beta \cdot \beta^{\sigma} = \epsilon$. Since $K$ is
imaginary, $\epsilon = 1$. Proposition \ref{Weil} applies to
ensure the existence of a polarized abelian surface $(A, \cL )$
defined over $\Q $ and isomorphic over $K$ to $(A_0, \cL_0)$.
Since $A\simeq A_0$ over $K$, we obtain that there is an
isomorphism  $\iota : \cO_D\stackrel{\sim }{\ra } \End_K(A)$.

Finally, we show that $\End^0_{\Q }(A)=\Q (\sqrt{m})$. The
triplets $(A, \iota , \cL )$ and $(A_0,\iota_0,\cL_0)$ are
isomorphic and the assertion $[A, \iota _{| R}, \cL ]\in X_D^{(m)}
(\Q )$ implies that for every $\sigma \in \Gal(\qbar /\Q )$ and
$\alpha \in R$, $\iota (\alpha )^{\sigma }=\iota (\alpha)\in \End
_{\qbar }(A)$. Hence $\iota (R)\subset \End_{\Q }(A)$. Thus $\Q
(\sqrt{m})=R\otimes \Q \subseteq \End^0_{\Q }(A)\subseteq
\End^0_{\qbar }(A)=B_D$. Since the only two subalgebras of $B_D$
that contain $\Q (\sqrt{m})$ are $\Q (\sqrt{m})$ and $B_D$
themselves, and the latter can not occur by Proposition
\ref{Thm0}, we obtain that $\End^0_{\Q }(A)=\Q (\sqrt{m})$, as we
claimed. $\Box $

\begin{remark} Let $(A, \cL )$ be a polarized abelian variety over
a number field $K$ and let $K_0\subseteq K$ be its field of
moduli. Assume that $\Aut (A, \cL )\simeq \{ \pm 1\}$. Then,
Proposition \ref{Weil} can be rephrased in the language of
cohomology of groups as follows. Let
$\mathrm{Br}_{K_0}=H^2(\Gal(\qbar /K_0), \qbar ^*)$ denote the
Brauer group of $K_0$. For any choice of isomorphisms $\{
\mu_{\sigma }, \sigma \in \Gal(\qbar /K_0)\}$ as above, the map
$$\begin{matrix}c: &\Gal(\qbar /K_0) \times \Gal(\qbar /K_0)&
\longrightarrow & \qbar ^* \\  & (\sigma , \tau ) & \mapsto &
\mu_{\sigma }\cdot \mu_{\tau}^{\sigma }\cdot \mu_{\sigma
\tau}^{-1}
\end{matrix}$$
produces a well-defined continuous cocycle $\xi (A, \cL )\in
\mathrm{Br}_{K_0}[2]$ in the $2$-torsion subgroup of
$\mathrm{Br}_{K_0}$. It is easily checked that different choices
of $\mu_{\sigma }$ for $\sigma \in \Gal(\qbar/K_0)$ lead to
cocycles that differ from a coboundary. Proposition \ref{Weil}
says that $(A, \cL )$ admits a model over $K_0$ if and only if
$\xi (A, \cL )=1\in \mathrm{Br}_{K_0}[2]$. By class field theory
(cf.\,\cite[Ch.\,X, \S 4, \S 5, \S 6]{Serre}), there is a natural
identification between $\mathrm{Br}_{K_0}[2]$ and the group
$(A_{K_0}, \otimes )$, where we let $A_{K_0}$ denote the set of
isomorphism classes of quaternion algebras over $K_0$. Upon this
identification, the proof of Theorem \ref{teodes} shows that if
$(A, \iota :\cO_D\stackrel{\sim }{\ra } \End _{\qbar }(A),\cL )$
is a polarized abelian surface over $\qbar $ with quaternionic
multiplication such that the field of moduli of $(A, \iota_{|R}
:R\hookrightarrow \End _{\qbar }(A),\cL )$ is $\Q $ for some order
$R\subset \cO_D$, $R\otimes \Q \simeq \Q (\sqrt{m})$ , then the
field of moduli of $(A, \iota , \cL )$ is an imaginary quadratic
field $K=\Q (\sqrt{-\delta })$ and $\xi (A, \cL ) = B_D\otimes
(\frac {-\delta , m}{\Q })$ in $\mathrm{Br}_{\Q }$.

\end{remark}

Motivated by the above result, we make the following definition.

\begin{defn}
The subset of descent points of $X_D^{(m)}(\Q )$ is $X_D^{(m)}(\Q
)_{d} :=$ $$\{ Q\in X_D^{(m)}(\Q )_{nh} : \pi_m^{-1}(Q)\subset
X_D(\Q (\sqrt{-\delta })), B_D\simeq (\frac{-\delta , m}{\Q })
\mbox{ for some }\delta >0\}.$$
\end{defn}
Set $r_m:=\# X_D^{(m)}(\Q)\,,\quad rh_m:=\#
X_D^{(m)}(\Q)_h\,,\quad rd_m:=  \# X_D^{(m)}(\Q)_{d}\,.$ $\\ $

Note that $r_m$ and $rd_m$ may be $+\infty$. The proof of Theorem
\ref{teodes} actually yields more information. Let $(D, m)$ be a
premodular pair over $\Q $. One may wonder how many abelian
surfaces $A/\Q $ exist up to isomorphism such that
$\End_{\Q}^0(A)=\Q(\sqrt m )$ and $\End_{\overline{\Q}}(A)=\cO_D$.
We make this precise in what follows.

\begin{defn}
Let $\cQ _{(D, m)}(\Q )$ denote the set of $\qbar $-isomorphism
classes of abelian surfaces $A$ defined over $\Q $ such that $\End
_{\qbar }(A)\simeq \cO _D$ and $\End _{\Q }^0(A)=\Q (\sqrt{m})$.
\end{defn}

As in Section 3, let $W_D$ denote the Atkin-Lehner group acting on
$X_D$. For a premodular pair $(D, m)$ over $\Q $, let
${W}_{D, m} = W_D/\langle \om _m\rangle $. As stated in
\cite[Proposition 3.2.2]{JoPh} or \cite[Proposition 5.5]{GoRo2},
the fixed points of an Atkin-Lehner involution $\om \in W_D$
acting on $X_D$ are Heegner points. Hence, the group $W_{D,m}$ is
naturally a subgroup of $\Aut _{\Q }(X_D^{(m)})$ which freely acts
on the set $X_D^{(m)}(\Q )_{d}$.

\begin{teor}\label{corr}

Let $(D, m)$ be a premodular pair over $\Q $. There is a
canonical one-to-one correspondence
$$
\cQ _{(D, m)}(\Q )\longleftrightarrow  {W}_{D,m}\backslash
X_D^{(m)}(\Q )_{d}
$$
and hence, if $D=p_1 \cdot ... \cdot p_{2r}$, then
$$|\cQ _{(D, m)}(\Q )|=\frac{rd_m}{2^{2r-1}}\,.$$
\end{teor}

{\em Proof.} Let $[A]\in \cQ _{(D, m)}(\Q )$ represented by an
abelian surface $A$ defined over $\Q $. In Part I of the proof of
Theorem \ref{teodes}, we show how can one attach a point $Q\in
X_D^{(m)}(\Q )_{d}$. Since, as explained in the last two
paragraphs of Section \ref{Shimura}, the group $W_{D,m}$ acts on
$Q =[A, \iota _{|R}, \cL ]$ by fixing the isomorphism class of $A$
and switching $\iota $ and $\cL $, we deduce that $A$ produces a
well-defined point in ${W}_{D,m}\backslash X_D^{(m)}(\Q )_d$. The
inverse map from ${W}_{D,m}\backslash X_D^{(m)}(\Q )_d$ onto $\cQ
_{(D, m)}(\Q )$ is constructed in Part II of the proof of Theorem
\ref{teodes} and for the same reason as above it does not depend
on the choice of $Q\in X_D^{(m)}(\Q )_d$ in its orbit under the
action by $W_{D,m}$. $\Box $

\begin{cor}\label{546}
Let $D>546$. Then there exist only finitely many $\qbar
$-isomorphism classes of abelian surfaces $A$ of $\GL_2$-type over
$\Q $ such that $\End_{\qbar }(A)\simeq \cO_D$.
\end{cor}

{\em Proof. } For a given discriminant $D$, assume that there
exist infinitely many abelian surfaces $A$ of $\GL_2$-type over
$\Q $ up to isomorphism over $\qbar $ such that $\End_{\qbar
}(A)$ $\simeq \cO_D$. According to Theorem \ref{MAIN} (i) and Theorem
\ref{corr}, there exists some $m\mid D$ such that $X_D^{(m)}$ has
infinitely many rational points over $\Q $. Since the degree of
the map $\pi_m:X_D\ra X_D^{(m)}$ is $2$, this implies that $X_D$
has infinitely many quadratic points. As is shown in \cite[Theorem
9]{Ro}, the largest such discriminant is $D=546$. $\Box $

\begin{remark}
One can check \cite[Table 3]{Ro} to find out which is precisely
the list of discriminants $D$ for which we can claim that
Corollary \ref{546} holds true.
\end{remark}

As a refinement of the above considerations, we wonder for which
pairs $(D, m)$ there exists a curve $C/ \Q $ of genus 2 such that
the Jacobian $J(C)$ has multiplication by $\Q(\sqrt{m})$ over $\Q
$ and quaternionic multiplication by $\cO_D$ over $\qbar $.

\begin{cor}\label{curve} Let $K=\Q(\sqrt{-\delta })$ be an imaginary quadratic
field. Then, there exists a curve $C/\Q $ defined over $\Q $ such
that $\End^0_{\Q }(J(C)) \simeq \Q (\sqrt{m})$ and $\End_{K}(J(C))$
$\simeq \cO_D$ if and only if $\pi_m(X_D(K))\cap X_D^{(m)}(\Q
)_{nh}\ne \emptyset $ and
$$B_D\simeq \left( \frac{-\delta,m}{\Q}\right)\simeq \left( \frac{-D,m}{\Q}\right)\,.$$
\end{cor}

{\em Proof.} By Theorem \ref{teodes}, there exists an abelian
surface $A/\Q $ such that $\End^0_{\Q }(A) \simeq \Q (\sqrt{m})$
and $\End_{K}(A) \simeq \cO_D$ if and only if $\pi_m(X_D(K))\cap
X_D^{(m)}(\Q )_{nh}\ne \emptyset $ and $B_D\simeq \left(
\frac{-\delta,m}{\Q}\right)\,.$ By Lemma \ref{pp}, $A$ admits a
principal polarization over $\Q $ if moreover $B_D\simeq \left(
\frac{-D,m}{\Q}\right)$. The result now follows because as is
well-known (cf.\,e.\,g.\,\cite[Theorem 3.1]{GoGuRo}), an
absolutely irreducible abelian surface $A/\Q $ is the Jacobian of
a smooth curve $C/\Q $ of genus $2$ if and only if $A$ is
principally polarizable over $\Q $. $\Box $

Next, we illustrate the above results with several examples.

\begin{example} In \cite[Lemma 4.5]{HaTs}, Hashimoto and Tsunogai provided
a family of curves of genus $2$ whose Jacobians have quaternionic
multiplication by $\cO_6$. These families specialize to infinitely
many curves defined over $\Q $. However, one can not expect that
to be always possible for a discriminant $D$ even when there is an
Atkin-Lehner quotient $X_D^{(m)}\simeq \PP ^1_{\Q }$. As we
pointed out in Section \ref{Conj}, computations due to Hasegawa
\cite{Ha} exhibit $B_{14}$ as a modular quaternion algebra. This
is indeed possible because $X_{14}^{(14)}\simeq \PP ^1_{\Q }$ but
there does not exist a curve $C/\Q$ of genus $2$ whose Jacobian
$J(C)$ is of $\GL _2$-type over $\Q $ and has quaternionic
multiplication by $\cO_{14}$ over $\qbar $, because $B_{14}\not
\simeq (\frac {-14, 2}{\Q })$, $(\frac {-14, 7}{\Q })$ nor $(\frac
{-14, 14}{\Q })$.
\end{example}

\begin{example}
An affine equation of the Shimura curve $X_6$ is $x^2+y^2+3=0$
(cf.\,Table 1) and the action of $w_6$ on this model is
$(x,y)\mapsto (-x,y)$. We have $X_6^{(6)}\simeq \PP ^1_{\Q }$ and
there exist infinitely many points on $X_6(K)$, $K=\Q
(\sqrt{-21})$, mapping to a rational point on $X_6^{(6)}$. Regard
$X_6$ as the coarse moduli space attached to $\hat {\cF }_{\cO_6,
\mu }$ for some $\mu\in \cO_6$ such that $\mu^2+6=0$. As stated in
\cite[Theorem 1.4.1]{JoPh}, the points $P\in X_6(K)$ are
represented by {\em principally} polarized abelian surfaces $P=[A,
\iota , \cL ]$ with quaternionic multiplication. Since $K$ splits
$B_6$, we can assume by \cite[Theorem 2.1.3]{JoPh} that $(A, \iota
, \cL )$ are defined over $K$. Moreover, since the class number
$h(K)>1$, we know by \cite[Theorem 3.1.5]{JoPh}, \cite[Theorem
5.11]{GoRo2} that there are no Heegner points on $X_6(K)$. As
explained in Section 2.1, this means that the abelian surfaces $A$
are absolutely irreducible. Finally, we have $\om_6(P) = \om_6[A,
\iota , \cL ] = [A, \mu^{-1}\iota \mu , \cL ]$ by \cite[Theorem
3.5]{Ro3} (cf.\,also \cite[Section 7, p.\,273]{Ro3}). Since
$\pi_6(P)\in X_6^{(6)}(\Q )$, this implies that $(A, \cL)$ are
isomorphic to their Galois conjugate $(A^{\sigma }, \cL ^{\sigma
})$, for $\Gal(K/\Q )=\langle \sigma \rangle $. However, since
$(\frac{-21,6}{\Q})\not\simeq B_6$, we conclude by Theorem
\ref{teodes} that there does not exist a model of $(A, \cL )$
defined over $\Q $.
\end{example}

\begin{example}
Let $f$ be the newform of  $S_2(\Gamma_0(243))$ with $q$-expansion
$$f=q+\sqrt 6 \,q^2+ 4\, q^4+\dots  \,.$$
The modular abelian surface $A_f$ obtained as an optimal quotient
of the Jacobian of $X_0(43)$ satisfies that $\End _{K}(A_f)$ is a
maximal order of the quaternion algebra $B_6$, where we let
$K=\Q(\sqrt{-3})$. By \cite[Theorem 7.1]{Ro1} we know that there
is a single class of principal polarizations $\cL_0$ on $A_f/K$ up
to $\qbar $-isomorphism. Hence $A_f\otimes K$ is the Jacobian of a
curve $C/K$. Since $\cL_0$ is isomorphic to its Galois conjugate
$\cL_0^{\sigma }$, it follows that $C_0$ is isomorphic to
$C_0^{\sigma }$ but, although its Jacobian $\Jac (C_0)=A_f\otimes
K$ admits a projective model over $\Q $, the curve $C_0$ can not
be defined over $\Q$ because $B_6\not\simeq (\frac{-6,6}{\Q})$. In
fact, by using similar methods to \cite{GoGuRo}, we obtain the
following equation for $C_0$:
\begin{center}
$y^2=(2+2\sqrt{-3})x^6+12 (-3+\sqrt{-3})x^5-12(3+7 \sqrt{-3})x^4+
4(69+7\sqrt{-3})x^3+ 8(-11+7 \sqrt{-3})x^2-18(1+5 \sqrt{-3})+12(2+\sqrt{-3}).$
\end{center}
It can be checked that its Igusa invariants are rational and there
is a morphism $\nu: C_0\rightarrow C_0^{\sigma}$ defined over $K$
such that $\nu^{\sigma}\cdot \nu$ is the hyperelliptic
involution.
\end{example}

\section{Rational points on quotient Shimura curves $X_D^{(m)}$ of
genus $\leq 1$}\label{g01}

The set of rational Heegner points $X_D^{(m)}(\Q)_h$ on an
Atkin-Lehner quotient of a Shimura curve is finite and its
cardinality $rh_m$ can be computed by using the following formula,
which stems from the work of Jordan on complex
multiplication\footnote{To be precise, the work in
\cite[Ch.\,3]{JoPh} restricts to complex multiplication by the full
ring of integers in imaginary quadratic fields. Since there exist
nonmaximal orders of class number $1$ or $2$, one actually needs to
apply the statements of \cite[Section 5]{GoRo2}.}.

\begin{prop}\label{CM} Let $D=p_1\cdot ...\cdot p_{2r}$ and let $m\vert D$. For $i=1$ or $2$,
let us denote by $\cR_i$ the set of orders of imaginary quadratic
fields whose class number is $i$. For any $R\in \cR_1$, set
$p_R=2$ when $\operatorname{disc}\, (R)=-2^k$ and $p_R$ to be the
single odd prime dividing $\operatorname{disc}(\,R)$, otherwise.
The number $rh_m$ is given by the following formula:
$$rh_m=\left\{
\begin{array}{{rc}}
2^{2r-1}\#\{ R\in\cR_1: \left( \frac{R}{p_i}\right)=-1\,\, \mbox{for all } p_i\vert D\}+{}\\
2^{2r-2}\#\{ R\in\cR_1: p_R\vert D\,, \left(
\frac{R}{p_i}\right)=-1\,\, \mbox{for all } p_i\vert
\frac{D}{p_R}\}+{}\\
\# \{R\in\cR_2: 2 D=-\operatorname{disc}\, (R)\}\phantom{{}+{}}&\mbox{, if $m=D$,}\\
\\
2^{r-2}\#\{ R\in\cR_1: p_R=\frac{D}{m}\}+\#\{R\in\cR_2:
-\frac{D}{\operatorname{disc\,}(R)}\in\Q^{*2}\}+{}\\
\frac{(-1)^p+1}{2} \#\{ R\in \cR_2:2D=-\operatorname{disc\,}(R)\}\phantom{{}+{}}&\mbox{, if $m=\frac{D}{p}$,}\\
\\
 0 &\mbox{otherwise.}
\end{array}
\right.
$$
where $(\frac{R}{p})=\begin{cases} (\frac{K}{p}) & \text{ if } p\nmid \text{cond}(R) \\
1 & \text{ if } p\mid \text{cond}(R) \end{cases}$ is the Eichler
symbol.

\end{prop}

{\em Proof. } This is a direct application of Proposition 5.5 and
Corollary 5.13 of \cite{GoRo2}. $\Box $

\vspace{0.2cm} In any case, $0\leq rd_m\leq r_m-rh_m$. The
condition $r_m=rh_m$ implies $rd_m=0$ and allows us to claim the
non existence of an abelian surface $A/\Q$ with $\End_{\Q}^0
A=\Q(\sqrt m)$ and $\End_{\overline{\Q}}^0 A=B_D$. Proving the
existence of such an abelian surface, i. e. $rd_m>0$, requires the
knowledge of an equation for $X_D$ and the action of $\omega_m$ on
it.

There are exactly twelve Shimura curves $X_D$ of genus $g \leq 2$.
For all of them,  $D=p\cdot q$ with $p,q $ primes and affine
equations for these curves are known  (cf.\ \cite{GoRo},
\cite{GoRo2}, \cite{Ku1}, \cite{Ku2}). For these equations, the
Atkin-Lehner involutions $\om_p,\om_q,\om_{p.q}$ act on the curve,
sending $(x,y)$ to $(-x,y)$, $(x,-y)$  and $(-x,y)$ in some
suitable order. The next table shows equations, genera and the
actions of $\om_p$ and $\om_q$ for these curves.

$$
\begin{array}{|c|c|c|c|c|c|}
\hline
D=p\cdot q& g&X_D & \om_p(x,y)&  \om_q(x,y) \\
\hline
2\cdot 3& 0& x^2+y^2+3=0 & (-x,-y)& (\phantom{-}x,-y)   \\
\hline
2\cdot 5& 0& x^2+y^2+2=0 & (\phantom{-}x,-y) &(-x,-y)   \\
\hline
2\cdot 11& 0&x^2+y^2+11=0 & (-x,-y)&(\phantom{-}x,-y)     \\
\hline
2\cdot 7& 1& (x^2-13)^2+7^3+2 y^2=0&  (-x,\phantom{-}y)&(-x,-y)    \\
\hline
3\cdot 5& 1& (x^2+3^5)(x^2+3)+3 y^2=0 & (-x,\phantom{-}y)&(-x,-y)   \\
\hline
3\cdot 7 & 1& x^4-658 x^2+7^6+7 y^2=0& (-x,-y)&(-x,\phantom{-}y)     \\
\hline
3\cdot 11& 1& x^4+30 x^2+3^8+3y^2=0 & (-x,\phantom{-}y)&(-x,-y)    \\
\hline
2\cdot 17& 1& 3x^4 -26  x^3+ 53 x^2+26 x+3+y^2=0 & (-\frac{1}{x},\frac{y}{x^2})&(-\frac{1}{x},\frac{-y}{x^2})    \\
\hline
2\cdot 23 & 1&(x^2-45)^2+23+2 y^2 = 0 &(-x,\phantom{-}y)&(-x,-y)   \\
\hline
2\cdot 13& 2&y^2 =-2 x^6+19 x^4- 24 x^2-169 &(-x,-y)&(-x,\phantom{-}y)   \\
\hline
2\cdot 19& 2&y^2=-16 x^6 -59 x^4-82 x^2-19   & (-x,-y)&(-x,\phantom{-}y)   \\
\hline
2\cdot 29& 2&2 y^2=  - x^6-39 x^4- 431x^2- 841   &(-x,-y)&(\phantom{-}x,-y)   \\
\hline
\end{array}
$$
\centerline{\ \ \ \ \ {\bf Table 1.} Equations and Atkin-Lehner
involutions on Shimura curves}
\medskip

\vskip 0.1 cm

Unfortunately, the construction of the above equations does not
allow us to distinguish the rational Heegner points among the
rational points on the curves $X_D^{(m)}$, unless they are fixed
by some Atkin-Lehner involution (however, for $D=6$ and $10$,
cf.\,\cite{Elk}). This forces the proof of next theorem to be more
elaborate.

\vskip 0.1 cm

\begin{teor}\label{triplets} For the twelve values of $D$ as above,
the  triplets $(r_m,rh_m,rd_m)$ take the following values:
$$\begin{array}{|c|c|c|c|}
\hline
D=p\cdot q &   (r_p,rh_p,rd_p)     &(r_q,rh_q,rd_q)   & (r_D,rh_D,rd_D) \\
\hline
2\cdot 3 & (\infty, 1, \infty)& (\infty,1,\infty)& (\infty,8,\infty)\\
\hline
2\cdot 5& (\infty,2, 0)   & (\infty,2, \infty) &(\infty,11,\infty)\\
\hline
2\cdot 7  &  (0, 0,0) & (6,2, 4)&(\infty,8,\infty )\\
\hline
2\cdot 11& (\infty,2, \infty)  &  (\infty,2, \infty) &(\infty,8,\infty)\\
\hline
2\cdot 13&(1,1,0) & (3,1,0)  &(\infty,10,>0) \\
\hline
2\cdot 17 &(0,0,0) & (0,0,0)  &(\infty,8,\infty) \\
\hline
2\cdot 19& (1,1,0)  & (3,1,0) &(\infty,8,>0)  \\
\hline
2\cdot 23 & (0, 0,0) & (2, 2,0)&  (\infty,8,\infty)\\
\hline
2\cdot 29 & (1,1,0)& (\infty,2, >0) &(\infty,13,>0)\\
\hline
3\cdot 5&  (\infty,2, 0 )   &  (4, 4,0) &(\infty,10, \infty)\\
\hline
3\cdot 7& (2, 2,0)  &(0, 0,0) &(\infty,10,\infty)\\
\hline
3\cdot 11&  (0,0,0) &(2,2,0)  &(\infty,10,\infty)\\
\hline
\end{array}
$$
\end{teor}
\centerline{\ \ \ \ \ {\bf Table 2.} Rational, Heegner and descent
points}
\medskip

\begin{proof}
We split the proof in three parts according to the genus $g$ of
$X_D$.

{\it Case $g=0$.} In all these cases the equation is $x^2+y^2=-d$,
for some prime $d\vert D$. For each pair $(D, m)$, the points on
$X_D(\sqrt{-\delta })$ of the form $(a, b\sqrt{-\delta})$, $(
b\sqrt{-\delta},a)$ or $(a\sqrt{-\delta}, b\sqrt{-\delta})$, with
$a,b\in\Q$ and a square-free integer $\delta \geq 1$, are the only
affine points on $X_D(\Q (\sqrt{-\delta }))$ which may provide
rational points on $X_D^{(m)}(\Q )$, depending on whether
$\omega_m$ maps $(x,y)$ to either $(x,-y)$, $(-x,y)$ or $(-x,-y)$.
Let $(\cdot , \cdot )$ denote the global Hilbert symbol over $\Q$.
When $\om _m(x, y)=(x, -y)$ or $(-x, y)$, such rational points
exist if and only if $(\delta,-d)=1$. Similarly, when $\om _m(x,
y)=(-x, -y)$, there exist points on $X_D(\Q (\sqrt{-\delta }))$
which project onto a rational point on $X_D^{(m)}(\Q )$ if and
only if $(-1,d\cdot \delta)=1$. It is easy to check that, for all
pairs $(D,m)\not =(10,2)$, these conditions and the descent
condition of Theorem \ref{teodes} have infinitely many solutions
for $\delta$. If we let $l$ be a prime number, we may take
$\delta$ as follows.
$$
\begin{array}{|c|c|c|}
\hline
(D,m) &  \delta &  \mbox{Conditions on }l\\
\hline
(6,2)& 3\,l &l\equiv 1\pmod{8}\\
\hline
(6,3), (6,6) & l & l\equiv 1\pmod{24}\\
\hline
(10,5) & 2\,l &l\equiv 1 \pmod {20}\\
\hline
(10,10) & 2\,l &l\equiv 1 \pmod {40}\\
\hline
(22,2)& 11\,l& l\equiv 1 \pmod 8\\
\hline
(22,11),(22,22)&l & l\equiv 1 \pmod {88}\\
\hline
\end{array}
$$
For the pair $(10,2)$, we have $d=2$ and the condition
$(\delta,-2)=1$ implies $5\not \vert \delta$ and thus
$(\frac{-\delta,2}{\Q})\not \simeq B_{10 }$, since
$(-\delta,2)_5=1$. Moreover, the two points $(1:\pm\sqrt{-1}:0)$
at infinity on the curve $X_{10}$ produce a rational point on
$X_{10}^{(2)}$ which is Heegner because its preimages are fixed
points by the Atkin-Lehner involution $\om_5$. We conclude that
for $D=10$, $rd_{2}=0$ despite $r_2=\infty$.

{\it Case $g=1$.} Let us first consider the cases for which $m\neq
D$. Note that $\omega_m$ does not act as $(x,y)\mapsto (-y,x)$.

Assume first that $D\neq 34$. Then, the genus of $X_D^{(m)}$ is
zero except for the cases in which $\omega_m$ maps $(x,y)$ to
$(-x,-y)$. The latter holds for the pairs $(D,m)=(14,7),(15,5),
(21,3)$, $(33,11)$ and $(46,23)$, and in these cases
$g(X_D^{(m)})=1$. To be more precise, the curves $X_D^{(m)}$ are
elliptic curves over $\Q$ and there is a single isogeny class of
conductor $D$ in each case. Their Mordell-Weil rank over $\Q $ is
$0$ and the orders of the group of rational torsion points on them
are $6$, $4$, $2$, $2$ and $2$, respectively. Only for the pair
$(D,m)=(14,7)$ do we have $r_m>rh_m$. But in this case, the two
rational Heegner points can be recognized from the affine equation
of the curve $X_{14}^{(7)}$ because they are fixed points by some
Atkin-Lehner involution. It then turns out that the four rational
non-Heegner points on $X_{14}^{(7)}$ are the projections of the
points $(\pm 8\sqrt{-1}, \pm 56\sqrt{-1}),(\pm 2\sqrt{-2},\pm 14
\sqrt{-2} )\in X_{14}(\overline{\Q})$. Since $\delta=1, 2$ satisfy
the descent condition, we have $rd_7=4$.

When $\omega_m$ acts as $(x,y)\mapsto (-x,y)$, we have $r_m=0$
except for $(D,m)=(15,3)$. The affine equation $(X+3^5)(X+3)+3
y^2=0$ for $X_{15}^{(3)}$ shows that there are no rational points
at infinity on this model. Moreover, it turns out that for all
$(X,y)\in X_{15}^{(3)}(\Q)$, the $5$-adic valuation of the
$X$-coordinate is $v_5(X)=0$. Since $\delta \equiv -X\pmod
{\Q^2}$, we have $5\not \vert \delta$. Thus, $(\frac{-\delta,
3}{\Q})\not\simeq B_{15}$ because $(-\delta, 3)_5=1$. It follows
that $rd_3=0$.

Finally, let us consider the particular case $D=34$. The equations
of the curves $X_{34}^{(2)}$ and $X_{34}^{(17)}$ are
$$
v^2+ 3 u^4-26 u^3+71 u^2-104 u+236=0\,, \quad  z^2+ 3 u^2-26
u+59=0\,,
$$
respectively, where $u=x-1/x$, $v=y(1+1/x^2)$ and $z=y/x$. Hence,
$X_{34}^{(2)}(\R )=\emptyset $ and $X_{34}^{(17)}(\R )=\emptyset $
and we obtain that $r_m=0$ in both cases. Note that this also
follows from Proposition~\ref{local}.

We now consider the case $m = D$.
We have $\omega_m:(x,y)\mapsto (x,-y)$ and the curve $X_D^{(D)}$ 
admits an affine model of the form $f(x)+ d Y=0$, where $f(x)\in \Q [x]$ 
is monic of degree $4$
and $(\frac{-d,D}{\Q})\simeq B_D$.
A point $(x_0, Y_0)\in
X_D^{(D)}(\Q)$ satisfies the descent condition if and only if
$(f(x_0),D)=1$, that is, $f(x_0)=u_0-D v_0^2 $ for some
$u_0,v_0\in\Q$. Hence, the descent condition for $(x_0, Y_0)$
turns out to be equivalent to the existence of a rational point on
the algebraic surface $S_D:f(x)=u^2-D v^2$ with $x=x_0$. For
$D=14,15,21,33,34$ and $46$ we have the following rational points
on $S_D$:
$(x_0,u_0,v_0)=(4,24,4),(1,44,8),(1,944,192),(8,145,16),( 9,  113,
18)$ and $(4,408,60)$, respectively. The elliptic curves $E_{D}:
f(x)=u^2-D v_0^2$ have at the least three rational points: two
rational points at infinity and the affine point $(x_0,u_0)$. It
can be easily checked that, for the above values of $D$, the
rational torsion subgroup $E_{D, \operatorname{tors}}(\Q )$ of
$E_{D}$ is of order $2$, except for~$D=34$, when $E_{D,
\operatorname{tors}}(\Q )$ has order~$1$. We conclude that the
Mordell-Weil rank of $E_D(\Q )$ is greater than $0$ and thus
$rd_D=\infty$.

{\it Case $g=2$.} The three curves $X_D$ are bielliptic. In
Cremona notation \cite[Table 1]{Cre}, the three elliptic quotients
$X_D/\langle \om_2\rangle$ are $26B2$, $38B2$, $58B2$, while
$X_{26}/\langle \om_{13}\rangle$, $X_{38}/\langle
\om_{19}\rangle$, $X_{58}/\langle \om_{58}\rangle$ are $26A1$,
$38A1$ and $58A1$, respectively.

For the three curves $X_D^{(2)}$, we have $rh_2=r_2$ and hence the
pairs $(D,2)$ are not premodular over $\Q$. For $(26,13)$ and
$(38,19)$, the single rational Heegner point corresponds to the
projection of the two points at  infinity because they are fixed
points by $\om_2$. The preimages of the other two rational points
are $(\pm\sqrt{-5},\pm 26)$ and $(\pm \sqrt{-5}/2,\pm 19/4)$ for
$D=26$ and $38$ respectively. In both cases $(-5,q)_5=-1$ and
hence $rd_q=0$. For the remaining cases, we have $r_m=\infty$ and
we can easily find $rh_m+1$ rational points on $X_D^{(m)}$
satisfying the descent condition.
\end{proof}

When $D$ is large enough so that the genus of $X_D$ is at least
$3$, we know no explicit equations describing $X_D$ and Theorem
\ref{teodes} can not be directly applied. For those pairs $(D,m)$
such that $X_D^{(m)}(\Q )\ne \emptyset $, one can still show that
$(D,m)$ is not premodular over $\Q $ provided $X_D^{(m)}(\Q
)$ is a finite set and one can prove that $X_D^{(m)}(\Q ) =
X_D^{(m)}(\Q )_{h}$. The cardinality of $X_D^{(m)}(\Q )_{h}$ can
be computed by Proposition \ref{CM}, whereas the computation of
$|X_D^{(m)}(\Q )|$ is a much more difficult problem.

There exist several Atkin-Lehner quotients $X_D^{(m)}$ which are
elliptic curves over $\Q $. It readily follows from comparing this
table with \cite[Table 1]{Cre} that six of them are elliptic curves
of rank zero (cf.\,\cite[Table 6.4]{RoPh} for the complete list
together with Weierstrass models for them)\footnote{Table 6.4 in
\cite{RoPh} contains a mistake: The list of elliptic curves
$X_D^{(m)}$ over $\Q $ tabled there is not complete, as
$X_{35}^{(7)}$, $X_{51}^{(3)}$ and $X_{115}^{(23)}$ are missing.
Indeed, these elliptic curves are $35A1$, $51A2$ and $115A1$,
respectively. This error was also reproduced in \cite[Table 2]{Ro},
where it was wrongly claimed that these curves failed to have
rational points over $\Q_5$, $\Q_{17}$ and $\Q_5$, respectively.}.
Namely, these are $X_{35}^{(7)}$, $X_{51}^{(3)}$, $X_{106}^{(53)}$,
$X_{115}^{(23)}$, $X_{118}^{(59)}$ and $X_{202}^{(101)}$, which
correspond to the elliptic curves $35A1$, $51A2$, $106D1$, $115A1$,
$118 D1$ and $202A_1$ respectively. In all cases but $(35, 7)$ the
number of rational points is $r_m=1$. Thus $rh_m=1$ since $r_m-rh_m$
is even whenever $r_m<\infty $. In particular, we get $rd_m=0$ and
thus $(51, 3)$, $(106,53)$, $(115, 23)$, $(118,59)$ and $(202,101)$
are not premodular over $\Q $. We can claim nothing about $(35,
7)$, since $r_m=3$, $rh_m=1$ and there is not an explicit equation
for $X_{35}$ at our disposal.

Combining the above with Theorems \ref{corr} and \ref{triplets},
together with those pairs exhibited in the last paragraph of Section
3, we obtain Theorem \ref{MAIN} (iii-iv-vi-vii) for all pairs but a
few ones which deserve more attention: namely, those $(D, m)$ such
that $X_D^{(m)}$ is a curve of genus $2$ and $X_D^{(m)}(\Q )\ne
\emptyset $. Computing the full list of rational points on these
curves is a harder task that we address in the next section.

\section{Covering techniques on bielliptic Shimura
curves $X_D^{(m)}$ of genus $2$}\label{g2}

It was shown in \cite[Proposition 4.2]{GoRo} that there exist
exactly ten Shimura curve quotients $X_D^{(m)}$ which are
bielliptic of genus $2$. Applying Proposition \ref{CM}, the
triplets $(D,m, rh_m)$ are $(91, 91,10)$, $(123, 123,10)$, $(141,
141,10)$, $(142, 2,0)$, $(142, 142,10)$, $(155, 155,10)$, $(158,
158,10)$, $(254, 254,8)$, $(326, 326,4)$ and $ (446, 446,6)$.

In this section we study the set of rational points on these
curves. We first list the $\Q$-rational points on each~$X_D^{(m)}$
which are easily found by a short search. Table~3 lists some small
rational points on the bielliptic curves $X_D^{m})$ of genus $2$.

\begin{table}
$$
\begin{array}{|c|c|cc|c|c|}
\hline D & m & &X_D^{(m)} & X_D^{(m)}(\Q)\\
\hline
91  & D & Y^2=&-X^6+19 X^4-3X^2+1  & (0,\pm 1), (\pm 1,\pm 4),(\pm 3,\pm 28)\\
\hline
123 & D & Y^2 =&   - 9X^6 + 19X^4+ 5 X^2+1&(0,\pm 1),(\pm 1,\pm 4),(\pm 1/3,\pm 4/3)\\
\hline
141  & D &  Y^2 =& 27X^6 - 5X^4 - 7X^2 + 1 & (\pm 1, \pm4), (\pm 1/3, \pm 4/9),\\
      &  &       &                         & (0, \pm 1), (\pm 11/13, \pm 4012/2197)\\
\hline
142  & 2 & Y^2 =& -16X^6 - 87X^4 - 146X^2 - 71 & \emptyset\\
\hline 142  & D & Y^2 =& 16X^6 + 9X^4 - 10X^2 + 1
     & \pm \infty ,(0, \pm 1), (\pm 1, \pm 4), (\pm 1/3, \pm 4/27)  \\
\hline 155  & D &  Y^2 = &25X^6 - 19X^4 + 11X^2 - 1
     &\pm \infty , (\pm1,\pm 4), (\pm 1/3,\pm 4/27)   \\
\hline 158  & D & Y^2 = &-8X^6 + 9X^4 + 14X^2 + 1
     & (\pm 1, \pm 4), (0, \pm 1), (\pm 1/3, \pm 44/27)  \\
\hline 254 & D &  Y^2 =& 8X^6 + 25X^4 - 18X^2 + 1
    & (0,\pm1),(\pm 1, \pm 2), (\pm 2, \pm 29) \\
\hline
326  & D &   Y^2 =& X^6+10 X^4 -63 X^2+4 &\pm \infty , (0,\pm 2) \\
\hline
446  & D &  Y^2 =& -16 X^6 - 7 X^4 + 38 X^2 + 1 &  (0,\pm 1),(\pm 1,\pm 4)\\
    \hline
\end{array}
$$
\centerline{\ \ \ \ \ {\bf Table 3.} Rational points on the
bielliptic $X_D^{(m)}$ of genus~$2$}
\end{table}

In this section, we show that Table~3 lists all rational points
for each~$X_D^{(m)}$. The case $D=142$, $m=2$ is straightforward:
there are no points in~$X_{142}^{(2)}({\mathbb R})$ from which it
follows that there are none in~$X_{142}^{(2)}({\mathbb Q})$. For
the other values of~$D,m$, each~$X_D^{(m)}$ has points everywhere
locally and so cannot be resolved in this way. We first recall the
techniques
from~\cite{bruin:283},\cite{bruin:thesis},\cite{fw1},\cite{flynn:qder},\cite{fw2},
which we summarize here in a simplified form adapted to the
curves~$X_D^{(m)}$. The fact that each~$X^{(m)}_D$ is bielliptic
allows a specialized version (cf.\,\cite{fw1}) of the same ideas
of \cite{bruin:crelle}; similar methods are available for
arbitrary hyperelliptic curves, as described
in~\cite{brufly:tow2cov} and~\cite{bruin:thesis}. Each of the
curves~$X_D^{(m)}$ is of genus~$2$ and of the form
$$
X_D^{(m)} :\ Y^2 = f_3 X^6 + f_2 X^4 + f_1 X^2 + f_0,\hbox{ with
}f_i\in \Z.
$$
Any such curve~$X_D^{(m)}$ has a map~$(X,Y) \mapsto (X^2,Y)$
from~$X_D^{(m)}$ to the elliptic curve~$Y^2 = f_3 w^3 + f_2 w^2 +
f_1 w + f_0$, and map~$(X,Y) \mapsto (1/X^2,Y/X^3)$
from~$X_D^{(m)}$ to the elliptic curve~$Z^2 = f_0 x^3 + f_1 x^2 +
f_2 x + f_3$. The Jacobian of~$X_D^{(m)}$ is~$\Q$-isogenous to the
product of these elliptic curves over~$\Q$ which, in all of these
examples, each have rank~$1$ (and no nontrivial torsion). It
follows that $\hbox{Jac}(X_D^{(m)})(\Q)$ has rank~$2$, and so
Chabauty techniques~\cite{chab:ratpoints} cannot be used, since
they only apply when the rank of the Mordell-Weil group of the
Jacobian is strictly less than the genus of the curve. It is
therefore necessary to imitate the technique in~\cite{fw1}, which
we briefly summarize here in a simplified form suited to these
examples. We first fix one of the above two elliptic curves -- it
does not matter which one; we shall use the latter elliptic curve,
since the resulting models will typically be slightly simpler.
Define~$E_D^{(m)}, (x_0, Z_0), \phi, t$ as follows.

\begin{align}
\label{eq:EDm}
 & E_D^{(m)} :\ Z^2 = f_0 x^3 + f_1 x^2 + f_2 x + f_3,\\
 \label{eq:x0Y0}
 & (x_0, Z_0) \hbox{ generates } E_D^{(m)}(\Q),\\
\label{eq:phi}
  & \phi :\ X_D^{(m)} \longrightarrow E_D^{(m)} :
        (X,Y) \mapsto (1/X^2,Y/X^3),\\
\label{eq:t}
  & t := \hbox{ root of } f_0 x^3 + f_1 x^2 + f_2 x + f_3,
\end{align}
so that $E_D^{(m)}(\Q)/2E_D^{(m)}(\Q) = \{ \infty, (x_0, Z_0) \}$.
Suppose that $(X,Y) \in X_D^{(m)}(\Q)$. Then, applying~$\phi$, we
let~$x = 1/X^2$ and $Z = Y/X^3$ so that $(x,Z) \in E_D^{(m)}(\Q)$.
We recall the injective homomorphism
$\mu:E_D^{(m)}(\Q)/2E_D^{(m)}(\Q) \rightarrow \Q(t)^*/(\Q(t)^*)^2$
defined by $\mu(\infty)=1$ and $\mu((x,Z))=f_0(x-t)$ from
\cite[Chapter X, Theorem 1.1]{Silverman}. It follows that
$\mu((x,Z))$ equals either~$1$ or~$f_0(x_0-t)$ in
$\Q(t)^*/(\Q(t)^*)^2$. Hence, either $f_0Z^2/(x-t)$ or
$(x_0-t)Z^2/(x-t)$ is a square. We can eliminate $Z^2$
using~(\ref{eq:EDm}) and simplification yields:
$$\begin{array}{rl}
\hbox{either }&f_0(f_0 x^2 + (f_0 t + f_1) x + (t^2 f_0 + t f_1 +
f_2))
              \in (\Q(t)^*)^2\\
\hbox{or }&(x_0 - t)(f_0 x^2 + (f_0 t + f_1) x + (t^2 f_0 + t f_1
+ f_2))
              \in (\Q(t)^*)^2.
\end{array}$$
Note that we do not really need $(x_0,Z_0)$; we only need the
square class of $x_0-t$. This can already be determined from the
$2$-Selmer group of $E_D^{(m)}$, without computing an explicit
generator of the Mordell-Weil group. In our examples, however, the
curve $E_D^{(m)}$ has small coefficients and finding an actual
generator is little more work.

Since~$x=1/X^2$ is a square itself, we can multiply either
quantity with~$x$ without changing its square class. Hence, we
have shown that if $(X,Y)\in X_D^{(m)}(\Q)$ then there exists
$y\in\Q(t)$ such that $(x,y)=(1/X^2,y)$ is a $\Q(t)$-rational
point on one of the curves
\begin{align}
\label{eq:FDm} F_D^{(m)} : y^2 &= f_0 x \bigl( f_0 x^2 + (f_0 t +
f_1) x + (t^2 f_0 + t f_1 + f_2)\bigr),
\\
\label{eq:GDm} G_D^{(m)} : y^2 &= (x_0-t) x \bigl(f_0 x^2 + (f_0 t
+ f_1)x + (t^2 f_0 + t f_1 + f_2)\bigr).
\end{align}
This gives a strategy for trying to prove that we have found all
of~$X_D^{(m)}(\Q)$; it is sufficient find all~$(x,y) \in
F_D^{(m)}(\Q(t))$ such that~$x\in\Q$ and similarly
for~$G_D^{(m)}$. This can be attempted using the techniques
in~\cite{bruin:283},\cite{bruin:thesis},\cite{fw1},\cite{flynn:qder},\cite{fw2},
which apply local techniques to bound the number of such points.
Since these articles already contain several worked examples of
this type, we merely provide a brief sketch of one case, to give
an idea of the general strategy, and to allow the reader to
interpret the tabular summary given later. This will be followed
by a description of any unusual features of difficult special
cases. The full details of the computations are available at:

\centerline{\begin{tabular}{ll}
\texttt{http://www.cecm.sfu.ca/{\textasciitilde}nbruin/shimura/}
\end{tabular}}

Consider~$X_{142}^{(142)} : Y^2 = 16X^6 + 9X^4 - 10X^2 + 1$, where
we wish to show that the only~$\Q$-rational points have
$X$-coordinate~$\infty, 0, 1, \pm 1/3$ (we use~$\infty$, depending
on context, as the notation both for the point at infinity and its
$X$-coordinate). Then~$E_{142}^{(142)}$ of~(\ref{eq:EDm}) is the
elliptic curve~$V^2 = x^3 - 10x^2 + 9x + 16$ over~$\Q$, which has
rank~1, with generator~$(x_0,V_0) = (1,4)$. Under~$(X,Y) \mapsto
(1/X^2, Y/X^3)$, the known points in~$X_{142}^{(142)}(\Q)$ map to
$(0,\pm 4)$, $\infty$, $(1,\pm 4)$ and $(9, \pm 4) = 2(1,\pm 4)$
in~$E_{142}^{(142)}(\Q)$. Letting~$t$ be the cubic number
satisfying~$t^3 - 10t^2 + 9t + 16$, the
curves~(\ref{eq:FDm}),(\ref{eq:GDm}) become
$$\begin{array}{rl}
F_{142}^{(142)} : y^2 &= x^3 + (t - 10)x^2 + (t^2 - 10t + 9)x,\\
G_{142}^{(142)} : y^2 &= (1-t)\bigl(x^3 + (t - 10)x^2 + (t^2 - 10t
+ 9)x\bigr).
\end{array}$$
Our known points in~$E_{142}^{(142)}(\Q)$ induce $\infty, (0,0),
(9,\pm (3 t^2 - 15 t)/4) \in F_{142}^{(142)}(\Q(t))$ and~$\infty,
(0,0), (1,\pm 4) \in G_{142}^{(142)}(\Q(t))$. It is sufficient to
show that there are no other points~$(x,y)$
in~$F_{142}^{(142)}(\Q(t))$ or~$G_{142}^{(142)}(\Q(t))$ for
which~$x\in \Q$. Note that, in each case~$\infty, (0,0)$ give the
entire torsion group, and so we have a point of infinite order.
Furthermore, a standard complete~$2$-descent or 2-isogeny descent,
as recently implemented by N.\ Bruin in Magma~\cite{magma} (or for
an older version, see~\cite{bruin:algae}), gives a Selmer bound
of~$1$ on the rank, and so both~$F_{142}^{(142)}(\Q(t))$
or~$G_{142}^{(142)}(\Q(t))$ have rank~$1$.

Consider, for example, the second curve~$G_{142}$. One can easily
check with finite field arguments, that~$R = (0,0) + 4(1,4) \in
G_{142}(\Q(t))$ is in the kernel of reduction modulo~$3$, which is
inert in~$\Q(t)/\Q$ (so that $\langle R \rangle$ is of finite
index in~$G_{142}(\Q(t))$), and can check that~$\langle R
\rangle$, $(0,0) + \langle R \rangle$, $(1,4) + \langle R \rangle$
and~$(1, -4) + \langle R \rangle$ and the only~4 cosets containing
possible~$(x,y) \in G_{142}^{(142)}(\Q(t))$ for which~$x\in \Q$
(note that, although~$(0,0),(1,4)$ do not
generate~$G_{142}(\Q(t))$, we do have that~$(0,0),(1-t,t^2-t-4)$
generate~$G_{142}(\Q(t))$ and that $(1,4) = (0,0) +
2(1-t,t^2-t-4)$, so that $(0,0),(1,4)$ can be treated as if they
are generators for the purposes of our $3$-adic argument). This
means that we only need to consider~$(x,y) = nR, (0,0) + nR, (1,4)
+ nR, (1,-4) + nR$, for~$n\in \Z_3$, and it is sufficient in each
case to show that~$n=0$ is the only case where~$x\in \Q$. Using
the formal exp and log functions in the~$3$-adic formal
group~\cite{fw1}, we can express $nR$ as~$\hbox{exp}(n\cdot
\hbox{log}(R))$ and deduce that
\medskip

\centerline{$1/(\hbox{$x$-coordinate of }n R) \equiv 8\cdot 3^2
n^2 + (2\cdot 3^2 n^2 + 3^3 n^4) t + (3^2 n^2 + 3^3 n^4)t^2 \hbox{
(mod $3^4$)}$,}
\medskip
\noindent where each coefficient of a power of~$t$ is a power
series in~$n$ defined over~$\Z_3$ whose coefficients converge
to~$0$. If~$x\in \Q$ then the coefficients of~$t$ and~$t^2$ must
be~$0$. Taking the coefficient of~$t$, we have a power series~$n$
for which~$n=0$ is a known double root, and for which the
coefficient of~$n^2$ has~$3$-adic absolute value strictly greater
than all subsequent coefficients of powers of~$n$. It follows
that~$n=0$ is the only solution. A~$3$-adic analysis of~$(0,0) +
nR, (1,4) + nR, (1,-4) + nR$ also shows that these can only have
$\Q$-rational $x$-coordinate when~$n=0$. We know that the
only~$(x,y) \in G_{142}^{(142)}(\Q(t))$ with~$x\in \Q$ are the
points with~$x=\infty, 0, 1$. Similarly, a~$5$-adic argument shows
that the only~$(x,y) \in F_{142}^{(142)}(\Q(t))$ with~$x\in \Q$
are the points with~$x=\infty, 0, 9$, as required.

These methods can be applied when the ranks of~$F_D^{(m)}(\Q(t))$
and~$G_D^{(m)}(\Q(t))$ are less than the degree of~$\Q(t)$, that
is, less than~$3$. Fortune is in our favour, since the ranks for
these examples indeed all turn out to be~$0,1$~or~$2$. Note that,
in the above rank~$1$ example, there was information to spare,
since either the coefficient of~$t$ or that of~$t^2$ could be used
to bound the number of solutions. For the rank~$2$ cases, one can
still obtain a bound, but the information from both power series
must be used.

The following table summarizes the computations.

$$
\begin{array}{|c|cc|c|c|c|}
\hline D = m & &F_D^{(m)}\hbox{ and }G_D^{(m)} & r & x \in \Q & p\\
\hline 91  & F_D^{(m)} : y^2=&x^3 + (t - 3)x^2 + (t^2 - 3t + 19)x
                          & 1 & \infty, 0, \frac{1}{9} & 5 \\
\hline 91  & G_D^{(m)} : y^2=&(1-t)\bigl( x^3 + (t - 3)x^2 + (t^2
- 3t + 19)x\bigr)
                          & 1 & \infty, 0, 1, 4 & 5 \\
\hline 123 & F_D^{(m)} : y^2=&x^3 + (t + 5)x^2 + (t^2 + 5t + 19)x
                           & 1 & \infty, 0, 9 & 5 \\
\hline 123 & G_D^{(m)} : y^2=&(1-t)\bigl(x^3 + (t + 5)x^2 + (t^2 +
5t + 19)x\bigr)
                           & 1 & \infty, 0, 1 & 7 \\
\hline 141 & F_D^{(m)} : y^2=&x^3 + (t - 7)x^2 + (t^2 - 7t - 5)x
                           & 1 & \infty, 0, 9 & 7 \\
\hline 141 & G_D^{(m)} : y^2=&(1-t)\bigl(x^3 + (t - 7)x^2 + (t^2 -
7t - 5)x\bigr)
                           & 2 & \infty, 0, 1, \frac{169}{121} & 7 \\
\hline 142 & F_D^{(m)} : y^2=&x^3 + (t - 10)x^2 + (t^2 - 10t + 9)x
                           & 1 & \infty, 0, 9 & 5 \\
\hline 142 & G_D^{(m)} : y^2=&(1-t)\bigl(x^3 + (t - 10)x^2 + (t^2
- 10t + 9)x\bigr)
                           & 1 & \infty, 0, 1 & 3 \\
\hline 155 & F_D^{(m)} : y^2=&x^3 + (t - 11)x^2 + (t^2 - 11t +
19)x
                           & 1 & \infty, 0, 9 & 7 \\
\hline 155 & G_D^{(m)} : y^2=&(t-1)\bigl(x^3 + (t - 11)x^2 + (t^2
- 11t + 19)x\bigr)
                           & 2 & \infty, 0, 1 & 3 \\
\hline 158 & F_D^{(m)} : y^2=&x^3 + (t + 14)x^2 + (t^2 + 14t + 9)x
                           & 2 & \infty, 0, 1, 9 & 5 \\
\hline 158 & G_D^{(m)} : y^2=&(-3-t)\bigl(x^3 + (t + 14)x^2 + (t^2
+ 14t + 9)x\bigr)
                           & 0 & \infty, 0 & \hbox{--} \\
\hline 254 & F_D^{(m)} : y^2=&x^3 + (t - 18)x^2 + (t^2 - 18t +
25)x
                           & 0 & \infty, 0 & \hbox{--} \\
\hline 254 & G_D^{(m)} : y^2=&(1-t)\bigl(x^3 + (t - 18)x^2 + (t^2
- 18t + 25)x\bigr)
                           & 2 & \infty, 0, 1, \frac{1}{4} & 5 \\
\hline 326 & F_D^{(m)} : y^2=& 4\bigl(4x^3 + (4t - 63)x^2 + (4t^2
- 63t + 10)x\bigr)
                           & 1 & \infty, 0 & 5 \\
\hline 326 & G_D^{(m)} : y^2=&-t\bigl(4x^3 + (4t - 63)x^2 + (4t^2
- 63t + 10)x\bigr)
                           & 0 & \infty, 0 & \hbox{--} \\
\hline 446 & F_D^{(m)} : y^2=&x^3 + (t + 38)x^2 + (t^2 + 38t - 7)x
                           & 1 & \infty, 0, 1 & 3 \\
\hline 446 & G_D^{(m)} : y^2=&(5-t)\bigl(x^3 + (t + 38)x^2 + (t^2
+ 38t - 7)x\bigr)
                           & 0 & \infty, 0 & \hbox{--} \\
\hline
\end{array}
$$
\centerline{\ \ \ \ \ {\bf Table 4.} Summary of computations}
\medskip

The second column gives the models for the curves~$F_D^{(m)}\hbox{
and }G_D^{(m)}$, and the third column gives the rank over~$\Q(t)$,
where the cubic number~$t$ is as defined in~(\ref{eq:t}). In all
cases, the torsion over~$\Q(t)$ consists only of~$\infty$
and~$(0,0)$. The fourth column gives the list of all~$x\in \Q$
which are~$x$-coordinates of a point~$(x,y)$ on the curve and
defined over~$\Q(t)$; the final column gives a prime~$p$ such that
a~$p$-adic argument proves that no other such~$x\in \Q$ are
possible. Of course, the rank~$0$ cases are trivial and so no such
prime is required.

The computations referenced above verify the following theorem.
\begin{teor}\label{thm:chabauty}
The curves $X_D^{(m)}$ have no $\Q$-rational points apart from
those given in Table 3.
\end{teor}

We conclude with mention of a few special features of the
computations. Recall that in the sketched worked example for the
case~$D=m=142$, it turned out that the Selmer bound from a
complete~$2$ descent was the same as the rank. However, for the
six cases~$F_{91}^{(91)}, G_{91}^{(91)}, G_{123}^{(123)},
F_{155}^{(155)}, F_{254}^{(254)}, G_{326}^{(326)}$, this bound is
two greater than the actual rank. In order to find a sharp bound,
one can perform a complete $2$-descent on the $2$-isogenous curve.
It follows that in each of these cases there are elements
of order~2 in the Shafarevich-Tate group over~$\Q(t)$.

In the other direction, there are two curves~$G_{155}^{(155)}$
and~$F_{326}^{(326)}$, where the group generated by images of the
known points in~$X_D^{(m)}(\Q)$ is less than the actual rank, so
one must search for the missing independent points of infinite
order. For example, the $2$-Selmer bound on the rank
of~$G_{155}^{(155)}(\Q(t))$ is~$2$, and the images of the known
points in~$X_{155}^{(155)}(\Q)$ give only~$\infty, (0,0), (1,4)$,
of which only~$(1,4)$ is of infinite order, so that we are missing
an independent point of infinite order. In this case, a naive
short search discovers the required point $((t^2 + 10t + 25)/4,
13t^2 - 11t + 10)$. The $2$-Selmer bound on the rank
of~$F_{326}^{(326)}(\Q(t))$ is~$1$, and the images of the known
points in~$X_{326}^{(326)}(\Q)$ do not give any points of infinite
order, and in fact the required point is

\centerline{ $(x,y) = \bigl( \frac{63540 t^2 - 1005167 t +
228495}{2888} ,
             \frac{90341332 t^2 - 1429154471 t + 325168047}{109744} \bigr).$}

This could not be found by a naive search, and we needed to use
the improved search techniques described in the appendix
of~\cite{bruin:thesis}, and recently implemented by N.~Bruin into
Magma~\cite{magma}.

As we explain in the last paragraph of Section 6,
Theorem~\ref{thm:chabauty} completes the proof of
Theorem~\ref{MAIN} (iii-iv-vi-vii). As we indicate at the end of
Section 3, parts (i-ii) of Theorem~\ref{MAIN} follow from Theorem
\ref{mod} and Proposition \ref{local}. Part (v) is proved in
Corollary \ref{546} and this gives the full proof of
Theorem~\ref{MAIN}.

From Theorem~\ref{thm:chabauty} and an analogue of Theorem
\ref{corr} we can also derive the following result.

\begin{cor}
For each of the pairs
$$(D, m) = (141, 141), (142, 142), (254, 254),$$
the number of $\qbar $-isomorphism classes of
abelian surfaces  $A/\qbar$ that admit an embedding $\iota :\Q
(\sqrt{m})\hookrightarrow \End _{\qbar }(A))$, whose field of
moduli is $\Q$, and such that $\End _{\qbar }(A)\simeq \cO _D$,
are $2$, $1$ and $1$, respectively.
\end{cor}


\begin{thebibliography}{99}

\bibitem{BoCa}
J.F.\ Boutot, H.\ Carayol, Uniformisation $p$-adique des courbes de
Shimura: les th\'eor\`emes de Cerednik et de Drinfeld. {\em
Ast\'erisque}, {\bf 196-197} (1991), 45-158.

\bibitem{bruin:283}
N.\ Bruin, The Diophantine equations $x^2 \pm y^4 = \pm z^6$ and
$x^2 + y^8 = z^3$, {\em Compositio Math.}, {\bf 118} (1999), no.\ 3,
305-321.

\bibitem{bruin:algae}
ALGAE, a program for 2-Selmer groups of elliptic curves over number
fields, available at
\texttt{http://www.cecm.sfu.ca/{\textasciitilde}bruin/ell.shar}

\bibitem{brufly:tow2cov}
N.\ Bruin and E.V.\ Flynn, Towers of $2$-covers of hyperelliptic
curves. Preprint PIMS-01-12 (2001), available at
\texttt{http://www.pims.math.ca/publications/\#preprints}

\bibitem{bruin:thesis}
N.\ Bruin, Chabauty methods and covering techniques applied to
generalized Fermat equations, Dissertation, University of Leiden,
Leiden, 1999. {\em CWI Tract}, {\bf 133}. Stichting Mathematisch
Centrum, Centrum voor Wiskunde en Informatica, Amsterdam (2002).

\bibitem{bruin:crelle}
N.\ Bruin, Chabauty methods using elliptic curves, {\em J. Reine
Angew. Math.}, {\bf 562} (2003), 27-49.

\bibitem{chab:ratpoints}
C.\ Chabauty, Sur les points rationnels des vari\'et\'es
alg\'ebriques dont l'irr\'egularit\'e est sup\'erieure \`a la
dimension, {\em C.\ R.\ Acad.\ Sci.\ Paris}, {\bf 212} (1941),
1022-1024.

\bibitem{Cl}
P. L.\ Clark, {\em Local and global points on moduli spaces of
abelian surfaces with potential quaternionic multiplication},
Harvard PhD.\ Thesis.


\bibitem{Cre} J.\ E.\ Cremona,  {\em Algorithms for modular
elliptic curves}, Cambridge Univ. Press, Cambridge, UK, 1992.


\bibitem{Cr}
J.\ E.\ Cremona, Abelian varieties with extra twist, cusp forms
and elliptic curves over imaginary quadratic fields, {\em J.\
London Math.\ Soc.\ } (2) ({\bf 45}) (1992), 401-416.

\bibitem{Elk} N.\ Elkies,  Shimura Curve Computations,
{\em Lect. Notes Comp.\  Sc.\ } {\bf 1423}, Proceedings of ANTS-3,
(1998); J.P.Buhler, ed.\ , 1-49.


\bibitem{DiRo}
L.\ Dieulefait, V.\ Rotger, The arithmetic of QM-abelian surfaces
through their Galois representations, {\em J.\ Algebra} {\bf 281}
(2004), 124-143.

\bibitem{DiRo2}
L.\ Dieulefait, V.\ Rotger, On abelian surfaces with potential
quaternionic multiplication, {\em Bull.\ Belg.\
Math.\ Soc.\ } {\bf 12:4}, (2005) 617-624.

\bibitem{Ed} B.\ Edixhoven, On the Andr\'e-Oort for Hilbert modular surfaces,
{\em Progress in Mathematics} {\bf 195} (2001), 133-155, Birkh\"auser.

\bibitem{fw1}
E.V.\ Flynn and J.L.\ Wetherell, Finding rational points on
bielliptic genus 2 curves, {\em Manuscripta Math.\ }, {\bf 100}
(1999), 519-533.

\bibitem{flynn:qder}
E.V.\ Flynn, On Q-derived polynomials, {\em Proc.\ Edinburgh Math.\
Soc.} {\bf 44} (2001), 103-110.

\bibitem{fw2}
E.V.\ Flynn and J.L.\ Wetherell, Covering collections and a
challenge problem of Serre. {\em Acta Arithmetica} {\bf XCVIII.2}
(2001), 197-205.

\bibitem{GLQ}
J.\ Gonz\'{a}lez, J.-C.\ Lario, J.\ Quer, Arithmetic of $\Q $-curves,
in {\em Modular curves and abelian varieties}, J.\ Cremona, J.-C.\
Lario, J.\ Quer, K.\ Ribet (eds.), Progress in Mathematics {\bf
224}, Birkh\"auser (2004), 125-140.


\bibitem{GoRo}
J.\ Gonz\'{a}lez, V.\ Rotger, Equations of Shimura curves of genus
two, {\em Intern.\ Math.\ Res.\ Not.\ } {\bf 14} (2004), 661-674.

\bibitem{GoRo2}
J.\ Gonz\'{a}lez, V.\ Rotger, Non-elliptic Shimura curves of genus
one, {\em J.\ Math.\ Soc.\ Japan} {\bf 58:4} (2006), 927-948.


\bibitem{GoGuRo}
J.\ Gonz\'{a}lez, J.\ Gu\`ardia, V.\ Rotger, Abelian surfaces of
$\GL_2$-type as Jacobians of curves, {\em Acta
Arithmetica} {\bf 116} (2005), 263-287.

\bibitem{Gr}
R.\ Greenberg, Seminar at Boston University and 
Letter to Frans Oort, unpublished.

\bibitem{Ha}
Y.\ Hasegawa, On some examples of modular QM-abelian surfaces, {\em
Proc.\ Japan Acad.\ , Ser.\ A } {\bf 72} (1996), 23-27.

\bibitem{HaTs}
K.\ Hashimoto, H. Tsunogai, On the Sato-Tate conjecture for
QM-curves of genus two, {\em Math.\ Comp.\ } {\bf 68} (1999),
1649-1662.

\bibitem{JoPh}
B.W.\ Jordan, {\em On the Diophantine arithmetic of Shimura curves},
Harvard PhD.\ Thesis (1981).


\bibitem{JoLi}
B.W.\ Jordan, R.\ Livn\'e, Local diophantine properties of Shimura
curves, {\em Math.\ Ann.\ } {\bf 270} (1985), 235-248.

\bibitem{Ku1}
A.\ Kurihara, On some examples of equations defining Shimura curves
and the Mumford uniformization, {\em J.\ Fac.\ Sci.\ Univ.\ Tokyo,
Sec.\ IA} {\bf 25} (1979), 277-301.

\bibitem{Ku2}
A.\ Kurihara, On $p$-adic Poincar\'{e} series and Shimura curves,
{\em Intern.\ J.\ Math.\ } {\bf 5} (1994), 747-763.

\bibitem{La}
S.\ Lang, {\em Algebraic Number Theory}, Graduate Texts in
Mathematics {\bf 110}, Springer (1970).

\bibitem{magma}
The Magma Computational Algebra System. Available from:

\newblock http://magma.maths.usyd.edu.au/magma/

\bibitem{Ma}
B.\ Mazur, Rational isogenies of prime degree, {\em Invent.\ Math.\
} {\bf 44} (1978), 129-162.


\bibitem{Mo}
B. J. J.\ Moonen, Models of Shimura varieties in mixed
characteristic, in {\em Galois represenations in arithmetic
geometry}, Eds.\ A. Scholl, R.\ Taylor, Cambridge University Press
(1998), 267-350.

\bibitem{Mu}
N.\ Murabayashi, On QM-abelian surfaces with a model of $\GL
_2$-type over $\Q $, Preprint.



\bibitem{Ogg}
A.P.\ Ogg, Mauvaise r\'eduction des courbes de Shimura, {\em
S\'{e}minaire de th\'{e}orie des nombres}, Progress in Mathematics
{\bf 59} Birkh\"{a}user Boston, Boston, MA (1983-84), 199-217.


\bibitem{Py}
E. E.\ Pyle, {\em Abelian varieties over $\Q $ with large
endomorphism algebras and their simple components over $\qbar $} in
{\em Modular curves and abelian varieties}, J.\ Cremona, J.-C.\
Lario, J.\ Quer, K.\ Ribet (eds.), Progress in Mathematics {\bf
224}, Birkh\"auser (2004), 189-239.

\bibitem{Ri}
K. A.\ Ribet, Abelian varieties over $\Q $ and modular forms, in
{\em Modular curves and abelian varieties}, J.\ Cremona, J.-C.\
Lario, J.\ Quer, K.\ Ribet (eds.), Progress in Mathematics {\bf
224}, Birkh\"auser (2004), 241-261.

\bibitem{RoPh}
V.\ Rotger, {\em On abelian varieties with quaternionic
multiplication and their moduli}, Universitat de Barcelona, PhD.\
Thesis (2003).

\bibitem{Ro}
V.\ Rotger, On the group of automorphisms of Shimura curves and
applications, {\em Compos.\ Math.\ } {\bf 132} (2002), 229-241.

\bibitem{Ro1}
V.\ Rotger, Quaternions, polarizations and class numbers, {\em J.\
Reine Angew.\ Math.\ } {\bf 561} (2003), 177-197.

\bibitem{Ro2}
V.\ Rotger, Modular Shimura varieties and forgetful maps, {\em
Trans.\ Amer. Math.\ Soc.\ } {\bf 356} (2004), 1535-1550.


\bibitem{Ro3}
V.\ Rotger, Shimura curves embedded in Igusa's threefold, in {\em
Modular curves and abelian varieties}, J.\ Cremona, J.-C.\ Lario,
J.\ Quer, K.\ Ribet (eds.), Progress in Mathematics {\bf 224},
Birkh\"auser (2004), 263-273.


\bibitem{Ro4}
V.\ Rotger, The field of moduli of quaternionic multiplication on
abelian varieties, {\em Intern.\ J.\ Math.\ Sc.\ } {\bf 52} (2004),
2795-2808.

\bibitem{RSY}
V.\ Rotger, A.\ Skorobogatov, A.\ Yafaev, Failure of the Hasse
principle for Atkin-Lehner quotients of Shimura curves over $\Q $,
{\em Moscow  Math.\ J.\ } {\bf 5:2}, (2005) 463-476.

\bibitem{Se}
J.-P.\ Serre, Sur les repr\'esentations modulaires de degr\'e $2$ de
$\Gal(\qbar /\Q )$, {\em Duke Math.\ J.\ } {\bf 54} (1987), 179-230.


\bibitem{Serre}
J.-P.\ Serre, {\em Corps Locaux}, Hermann, Paris (1968).

\bibitem{Sh67}
G.\ Shimura, Construction of class fields and zeta functions of
algebraic curves, {\em Ann.\ Math.\ } {\bf 85} (1967), 58-159.

\bibitem{Sh2}
G.\ Shimura, On the real points of an arithmetic quotient of a
bounded symmetric domain, {\em Math.\ Ann.\ } {\bf 215} (1975),
135-164.

\bibitem{Silverman}
J.H.\ Silverman, {\it The arithmetic of elliptic curves}, Graduate
Texts in Mathematics {\bf 106}, Springer, New York~(1986).

\bibitem{Vi}
M.F.\ Vign\'{e}ras, {\em Arithm\'{e}tique des alg\`{e}bres de
quaternions}, Lecture Notes in Mathematics {\bf 800}, Springer
(1980).

\bibitem{We}
A.\ Weil, The field of definition of a variety, {\em Amer.\ J.\
Math.\ } {\bf 78} (1956), 509-524.

\end{thebibliography}
\end{document}